\documentclass[11pt]{article}
\usepackage[english]{babel}
\usepackage{amsmath,amsthm,amssymb}
\usepackage{mathrsfs}
\usepackage{bbm}
\usepackage{amsfonts}
\usepackage[margin=0.8in]{geometry}
\usepackage{graphicx}
\usepackage{caption}
\usepackage{subcaption}
\usepackage{wasysym}
\usepackage{enumerate}
\usepackage{algorithm2e}
\usepackage{tabularx}
\usepackage{color}
\usepackage{wrapfig}
\usepackage{todonotes}

\newtheorem{thm}{Theorem}[section]
\newtheorem{ass}[thm]{Assumption}
\newtheorem{cor}[thm]{Corollary}
\newtheorem{lem}[thm]{Lemma}
\newtheorem{prop}[thm]{Proposition}
\newtheorem*{hyp*}{Hypothesis}
\theoremstyle{definition}
\newtheorem{defn}[thm]{Definition}

\theoremstyle{rem}
\newtheorem{rem}[thm]{Remark}

\numberwithin{equation}{section}
\newcommand{\R}{\mathbb R}

\newcommand{\bbF}{\mathbb F}

\newcommand{\mcA}{\mathcal{A}}

\newcommand{\mcB}{\mathcal{B}}
\newcommand{\mcD}{\mathcal D}
\newcommand{\mcL}{\mathcal L}

\newcommand{\mcF}{\mathcal F}

\newcommand{\mcT}{\mathcal T}
\newcommand{\mcP}{\mathcal P}
\newcommand{\mcO}{\mathcal O}
\newcommand{\mcM}{\mathcal M}
\newcommand{\mcH}{\mathcal H}

\newcommand{\mcU}{\mathcal U}
\newcommand{\mcS}{\mathcal S}
\newcommand{\mcY}{\mathcal Y}
\newcommand{\mcZ}{\mathcal Z}

\newcommand{\E}{\mathbb{E}}
\newcommand{\Prob}{\mathbb{P}}

\newcommand{\esssup}{\mathop{\rm{ess}\sup}}
\newcommand{\essinf}{\mathop{\rm{ess}\inf}}

\newcommand{\ett}{\mathbbm{1}}
\newcommand{\Vlow}{V_-}
\newcommand{\Vup}{V_+}
\newcommand{\Vlowh}{V_{-,h}}
\newcommand{\Vuph}{V_{+,h}}

\newcommand{\caglad}{c\`agl\`ad~}

\newcommand{\trace}{{\rm Tr}}

\newcommand{\ie}{\textit{i.e.\ }}
\newcommand{\eg}{\textit{e.g.\ }}
\newcommand{\etal}{\textit{et.~al.\ }}

\begin{document}

\title{Zero-sum Stochastic Differential Games of Impulse Versus Continuous Control by FBSDEs\footnote{This work was supported by the Swedish Energy Agency through grant number 48405-1}}

\author{Magnus Perninge\footnote{M.\ Perninge is with the Department of Physics and Electrical Engineering, Linnaeus University, V\"axj\"o,
Sweden. e-mail: magnus.perninge@lnu.se.}} %
\maketitle
% ----------------------------------------------------------------
\begin{abstract}
We consider a stochastic differential game in the context of forward-backward stochastic differential equations, where one player implements an impulse control while the opponent controls the system continuously. Utilizing the notion of ``backward semigroups'' we first prove the dynamic programming principle (DPP) for a truncated version of the problem in a straightforward manner. Relying on a uniform convergence argument then enables us to show the DPP for the general setting. In particular, this avoids technical constraints imposed in previous works dealing with the same problem. Moreover, our approach allows us to consider impulse costs that depend on the present value of the state process in addition to unbounded coefficients.

Using the dynamic programming principle we deduce that the upper and lower value functions are both solutions (in viscosity sense) to the same Hamilton-Jacobi-Bellman-Isaacs obstacle problem. By showing uniqueness of solutions to this partial differential inequality we conclude that the game has a value.
\end{abstract}

% ----------------------------------------------------------------
\section{Introduction}

The history of differential games is almost as long as the history of modern optimal control theory and traces back to the seminal work by Isaacs~\cite{Isaacs65}. To counter the unrealistic idea that one of the players have to give up their control to the opponent, Elliot and Kalton introduced the notion of strategies defined as non-anticipating maps from the opponents set of controls to the players own controls~\cite{ElliotKalton72}. Assuming that one player plays a strategy while the opponent plays a classical control, Evans and Souganidis \cite{Evans84} used the theory of viscosity solutions to find a representation of the upper and lower value functions in deterministic differential games as solutions to Hamilton-Jacobi-Bellman-Isaacs (HJBI) equations. Using a discrete time approximation technique, this was later translated to the stochastic setting by Flemming and Souganidis \cite{FlemSoug89}. The natural terminology for these games being zero-sum stochastic differential games (SDGs). Using the theory of backward stochastic differential equations (BSDEs), in particular the notion of backward semigroups, Buckdahn and Li~\cite{Buckdahn08} simplified the arguments and further extended the results in \cite{FlemSoug89} to cost functionals defined in terms of BSDEs.

Just as stochastic control was extended to various types of controls in the latter half of the previous century (notably to controls of impulse type in~\cite{BensLionsImpulse}), so has stochastic differential games. Tang and Hou~\cite{TangHouSWgame} considered the setting of two-player, zero-sum SDGs where both players play switching controls (a particular type of impulse control). Their result was later extended by Djehiche \etal \cite{BollanSWG1,DjehicheSWG2} to incorporate stochastic switching-costs. In the context of general impulse controls, Cosso~\cite{Cosso13} considered a zero-sum game where both players play impulse controls. By adapting the theory developed in~\cite{Buckdahn08}, L.~Zhang recently extended these results to cost functionals defined by BSDEs~\cite{LZhang21}.

%(sometimes also referred to as a \emph{stochastic})

In the present work we will be dealing with SDGs where one player plays an impulse control while the opponent plays a continuous control. This type of game problems have previously be considered by Azimzadeh~\cite{Azimzadeh19} for linear expectations and when the intervention costs are deterministic and by Bayraktar \etal \cite{BayraktarRobust} when the impulse control is of switching type. We follow the path described above where the cost functional is defined in terms of the solution to a BSDE and introduce the lower value function
\begin{align*}
  \Vlow(t,x):=\essinf_{\alpha^S\in\mcA^S_t}\esssup_{u\in\mcU_t}J(t,x;u,\alpha^S(u))
\end{align*}
and the upper value function
\begin{align*}
  \Vup(t,x):=\esssup_{u^S\in\mcU^S_t}\essinf_{\alpha\in\mcA_t}J(t,x;u^S(\alpha),\alpha)
\end{align*}
with $J(t,x;u,\alpha):=Y^{t,x;u,\alpha}_t$, where the pair $(Y^{t,x;u,\alpha},Z^{t,x;u,\alpha})$ solves the non-standard BSDE
\begin{align}\nonumber
Y^{t,x;u,\alpha}_s&=\psi(X^{t,x;u,\alpha}_T)+\int_s^Tf(r,X^{t,x;u,\alpha}_r,Y^{t,x;u,\alpha}_r,Z^{t,x;u,\alpha}_r,\alpha_r)dr
\\
&\quad-\int_s^T Z^{t,x;u,\alpha}_rdW_r - \Xi^{t,x;u,\alpha}_{T+}+\Xi^{t,x;u,\alpha}_s. \label{ekv:non-ref-bsde}
\end{align}
In the above definitions, $\mcU$ (resp. $\mcA$) and $\mcU^S$ (resp. $\mcA^S$) represent the set of impulse (resp. continuous) controls and their corresponding non-anticipative strategies. The generic member of $\mcU$ will be denoted by $u:=(\tau_i,\beta_i)_{1\leq i\leq N}$ where $\tau_i$ is the time of the $i^{\rm th}$ intervention and $\beta_i$ is the corresponding impulse, taking values in the compact set $U$. Moreover, the impulse cost process $\Xi$ is defined as
\begin{align}
\Xi^{t,x;u,\alpha}_s:=\sum_{j=1}^N\ett_{[\tau_j<s]}\ell(\tau_j,X^{t,x;[u]_{j-1},\alpha}_{\tau_{j}},\beta_j),
\end{align}
%Thus,
%\begin{align}
%\Xi^{t,x;u,\alpha}_{T+}-\Xi^{t,x;u,\alpha}_s=\sum_{j=1}^N\ett_{[\tau_j\geq s]}\ell(\tau_j,X^{t,x;[u]_{j-1},\alpha}_{\tau_{j}},\beta_j).
%\end{align}
where $[u]_j:=(\tau_i,\beta_i)_{1\leq i\leq N\wedge j}$ and $X^{t,x;u,\alpha}$ solves the impulsively and continuously controlled SDE
\begin{align}
X^{t,x;u,\alpha}_s&=x+\int_t^s a(r,X^{t,x;u,\alpha}_r,\alpha_r)dr+\int_t^s\sigma(r,X^{t,x;u,\alpha}_r,\alpha_r)dW_r\label{ekv:forward-sde1}
\end{align}
for $s\in [t,\tau_{1})$ and
\begin{align}
X^{t,x;u,\alpha}_{s}&=\Gamma(\tau_j\vee t,X^{t,x;[u]_{j-1},\alpha}_{\tau_j\vee t},\beta_j)+\int_{\tau_j\vee t}^s a(r,X^{t,x;u,\alpha}_r,\alpha_r)dr+\int_{\tau_j\vee t}^s\sigma(r,X^{t,x;u,\alpha}_r,\alpha_r)dW_r,\label{ekv:forward-sde2}
\end{align}
whenever $s\in [\tau_{j},\tau_{j+1})$ with $\tau_{N+1}:=\infty$.

We show that $\Vlow$ and $\Vup$ are both viscosity solutions to the Hamilton-Jacobi-Bellman-Isaacs quasi-variational inequality (HJBI-QVI)
\begin{align}\label{ekv:var-ineq}
\begin{cases}
  \min\{v(t,x)-\mcM v(t,x),-v_t(t,x)-\inf_{\alpha\in A}H(t,x,v(t,x),Dv(t,x),D^2v(t,x),\alpha)\}=0,\\
  \quad\forall (t,x)\in[0,T)\times \R^d \\
  v(T,x)=\psi(x),
\end{cases}
\end{align}
where $\mcM v(t,x):=\sup_{b\in U}\{v(t,\Gamma(t,x,b))-\ell(t,x,b)\}$ and
\begin{align*}
  H(t,x,y,p,X,\alpha):=p\cdot a(t,x,\alpha)+\frac{1}{2}\trace [\sigma\sigma^\top(t,x,\alpha)X]+f(t,x,y,p^\top \sigma(t,x,\alpha),\alpha).
\end{align*}
We then move on to prove that \eqref{ekv:var-ineq} admits at most one solution, leading to the main contribution of the paper, namely the conclusion that the game has a value, \ie that $\Vlow\equiv \Vup$.

As in most previous works on stochastic differential games involving impulse controls, the main technical difficulty we face is showing continuity of the upper and lower value functions in the time variable. In previous works such as \cite{TangHouSWgame,Cosso13,FZhang11} continuity is simplified by assuming that the intervention costs do not depend on the state and are non-increasing in time. In \cite{Azimzadeh19} the assumption of non-increasing intervention costs is replaced by one where the impulse player commits to, at the start of the game, limit to a fixed number of $q\geq 0$ impulses (where $q$ can be chosen arbitrarily large) in addition to assuming that impulses can only be made at rational times.

In the present work we take a completely different approach to the above mentioned articles, where we first show continuity under a truncation and then show that the truncated value functions converge uniformly to the true value functions on compact sets.

The paper is organized as follows. In the next section we give some preliminary definitions and describe the by now well established theory of viscosity solutions to partial differential equations (PDEs) as well as the notion of backward semigroups. Then, in Section~\ref{sec:FBSDEs} we give some preliminary estimates on the solutions to the non-standard BSDE in \eqref{ekv:non-ref-bsde}. Section~\ref{sec:DPP} is devoted to showing that dynamic programming principles hold for the lower and upper value functions. The proof that the lower and upper value functions are both solutions in viscosity sense to the same HJBI-QVI, that is \eqref{ekv:var-ineq}, is given in Section~\ref{sec:hjbi-qvi} while the uniqueness proof is postponed to Section~\ref{sec:unique}.

%\fbox{\parbox{0.9\textwidth}{
%The task is to:
%\begin{enumerate}
%  \item Derive the dynamic programming relations
%  \begin{align*}
%  \Vlow(t,x)=\essinf_{\alpha^S\in\mcA^S_{t,t+h}}\esssup_{u\in\mcU_{t,t+h}}G_{t,t+h}^{t,x;u,\alpha^S(u)}[\Vlow(t+h,X^{t,x;u,\alpha^S(u)}_{t+h})]
%  \end{align*}
%  and
%  \begin{align*}
%  \Vup(t,x)=\esssup_{u^S\in\mcU^S_{t,t+h}}\essinf_{\alpha\in\mcA_{t,t+h}}G_{t,t+h}^{t,x;u^S(\alpha),\alpha}[\Vup(t+h,X^{t,x;u,\alpha^S(u)}_{t+h})],
%  \end{align*}
%  for $t\in[0,T)$ and $h\in[0,T-t]$.
%  \item Show that \eqref{ekv:var-ineq} has at most one viscosity solution.
%  \item Show that $\Vlow$ and $\Vup$ are both viscosity solutions to \eqref{ekv:var-ineq}.
%  \item Show that $\Vlow(t,x)=\esssup_{\alpha^S\in\mcA^S}\bar Y^{t,x,\alpha^S(u^*)}_t$.
%\end{enumerate}
%}
%}
%\bigskip

\section{Preliminaries\label{sec:prel}}
We let $(\Omega,\mcF,\Prob)$ be a complete probability space on which lives a $d$-dimensional Brownian motion $W$. We denote by $\bbF:=(\mcF_t)_{0\leq t\leq T}$ the augmented natural filtration of $W$.\\

\noindent Throughout, we will use the following notation:
\begin{itemize}
  \item $\mcP_{\bbF}$ is the $\sigma$-algebra of $\bbF$-progressively measurable subsets of $[0,T]\times \Omega$.
  \item For $p\geq 1$, we let $\mcS^{p}$ be the set of all $\R$-valued, $\mcP_{\bbF}$-measurable \caglad processes $(Z_t: t\in [0,T])$ such that $\|Z\|_{\mcS^p}:=\E\big[\sup_{t\in[0,T]} |Z_t|^p\big]<\infty$ and we let $\mcS^p_c$ be the subset of processes that are continuous.
  \item We let $\mcH^{p}$ denote the set of all $\R^d$-valued $\mcP_{\bbF}$-measurable processes $(Z_t: t\in[0,T])$ such that $\|Z\|_{\mcH^p}:=\E\big[\big(\int_0^T |Z_t|^2 dt\big)^{p/2}\big]^{1/p}<\infty$.
  \item We let $\mcT$ be the set of all $\bbF$-stopping times and for each $\eta\in\mcT$ we let $\mcT_\eta$ be the corresponding subsets of stopping times $\tau$ such that $\tau\geq \eta$, $\Prob$-a.s.
  \item We let $\mcA$ be the set of all $A$-valued processes $\alpha\in \mcH^2$ where $A$ is a compact set.
  %\item For each $\tau\in\mcT$, we let $\mcA(\tau)$ be the set of all $\mcF_\tau$-measurable random variables taking values in $U$.
  %\item For any $k\geq 0$, we let $\bar\mcT^k:=\{(\tau_1,\ldots,\tau_k)\in\mcT^k:\tau_1\leq \cdots\leq\tau_k\}$.
  \item We let $\mcU$ be the set of all $u=(\tau_j,\beta_j)_{1\leq j\leq N}$, where $(\tau_j)_{j=1}^N$ is a non-decreasing sequence of $\bbF$-stopping times and $\beta_j$ is a $\mcF_{\tau_j}$-measurable r.v.~taking values in $U$, such that $\Xi^{t,x;u,\alpha}_T\in L^2(\Omega,\mcF_T,\Prob)$ for all $\alpha\in\mcA$.
  %\item We let $\mcU^f$ denote the subset of $u\in\mcU$ for which $N(t):=\sup\{j:\tau_j\leq t\}$ is $\Prob$-a.s.~finite on compacts (\ie $\mcU^f:=\{u\in\mcU:\: \Prob\left[\{\omega\in\Omega : N(t)>k, \:\forall k>0\}\right]=0, \:\forall t\in\R_+\}$) and for all $k\geq 0$ we let $\mcU^k:=\{u\in\mcU:\:N\leq k,\,\Prob{\rm - a.s.}\}$.
  \item For stopping times $\underline\eta\leq \bar\eta$ we let $\mcU_{\underline\eta,\bar\eta}$ be the subset of $\mcU$ with $\underline\eta\leq\tau_j\leq\bar\eta$, $\Prob$-a.s.~for $j=1,\ldots,N$. Similarly, we let $\mcA_{\underline\eta,\bar\eta}$, be the restriction of $\mcA$ to all $\alpha:\Omega\times [\underline\eta,\bar\eta]\to A$. When $\bar\eta=T$ we use the shorthands $\mcU_{\underline\eta}$ and $\mcA_{\underline\eta}$.
  \item For any $u\in\mcU$, we let $[u]_{j}:=(\tau_i,\beta_i)_{1\leq i\leq N\wedge j}$. Moreover, we introduce $N(s):=\max\{j\geq 0:\tau_j\leq s\}$ and let $u_s:=[u]_{N(s)}$ and $u^s:=(\tau_j,\beta_j)_{N(s)+1\leq j\leq N}$.
  \item We let $\Pi_{pg}$ denote the set of all functions $\varphi:[0,T]\times\R^n\to\R$ that are of polynomial growth in $x$, \ie there are constants $C,\rho>0$ such that $|\varphi(t,x)|\leq C(1+|x|^\rho)$ for all $(t,x)\in [0,T]\times\R^n$.
\end{itemize}

We also mention that, unless otherwise specified, all inequalities between random variables are to be interpreted in the $\Prob$-a.s.~sense.

\begin{defn}
We introduce that notion of \emph{non-anticipative strategies} defined as all maps $u^S:\mcA\to\mcU$ for which $(u^S(\alpha))_s=(u^S(\tilde \alpha))_s$ whenever $\alpha_r=\tilde\alpha_r$, $d\lambda\times d\Prob$-a.e. on $[0,s]\times\Omega$ (resp. $\alpha^S:\mcU\to\mcA$ for which $(\alpha^S(u))_s=(\alpha^S(\tilde u))_s$ whenever $\tilde u_s=u_s$, $\Prob$-a.s.). We denote by $\mcU^S$ (resp. $\mcA^S$) the set of non-anticipative strategies.

Moreover, we define the restrictions to an interval $[\underline\eta,\bar\eta]$ denoted $\mcU^S_{\underline\eta,\bar\eta}$ (resp. $\mcA^S_{\underline\eta,\bar\eta}$) as all non-anticipative maps $u^S:\mcA_{\underline\eta,\bar\eta}\to\mcU_{\underline\eta,\bar\eta}$ (resp. $\alpha^S:\mcU_{\underline\eta,\bar\eta}\to\mcA_{\underline\eta,\bar\eta}$).
\end{defn}

\begin{defn}
We will rely heavily on approximation schemes where we limit the number of interventions in the impulse control. To this extent we let $\mcU^k:=\{u\in\mcU:N\leq k,\,\Prob-{\rm a.s.}\}$ for $k\geq 0$ and let $\mcU^{S,k}$ be the corresponding set of non-anticipative strategies $u^S:\mcA\to \mcU^k$.
\end{defn}

\begin{defn}
We introduce the concatenation of impulse controls $\oplus$ as
\begin{align*}
  (\tau_j,\beta_j)_{1\leq j\leq N}\oplus (\tilde\tau_j,\tilde\beta_j)_{1\leq j\leq \tilde N}:=((\tau_1,\beta_1),\ldots,(\tau_N,\beta_N),(\tilde \tau_1\vee\tau_N,\tilde\beta_1),\ldots,(\tilde\tau_{\tilde N}\vee\tau_N,\beta_N))
\end{align*}
and note that for each $\eta\in \mcT$ we have the decomposition $u=u_\eta\oplus u^\eta$.

Similarly, when $0\leq t\leq s\leq T$ we let the concatenation of $\alpha\in\mcA_{t,s}$ and $\tilde\alpha\in \mcA_{s}$ at $s$ be defined as
\begin{align*}
  (\alpha\oplus_{s}\tilde\alpha)_r:=\ett_{[t,s)}(r)\alpha_r+\ett_{[s,T]}(r)\tilde\alpha_r
\end{align*}
for all $r\in [t,T]$.
\end{defn}

Throughout, we make the following assumptions on the parameters in the cost functional where $C>0$ and $\rho>0$ are fixed constants:
\begin{ass}\label{ass:on-coeff}
\begin{enumerate}[i)]
  \item\label{ass:on-coeff-f} We assume that $f:[0,T]\times \R^n\times\R\times\R^{d}\times A$ is Borel measurable, of polynomial growth in $x$, \ie there is a $C>0$ and a $\rho\geq 0$ such that
  \begin{align*}
    |f(t,x,0,0,\alpha)|\leq C(1+|x|^\rho)
  \end{align*}
  for all $\alpha\in A$, and that there is a constant $k_f>0$ such that for any $t\in[0,T]$, $x,x'\in\R^n$, $y,y'\in\R$, $z,z'\in\R^{d}$ and $\alpha\in A$ we have
  \begin{align*}
    |f(t,x',y',z',\alpha)-f(t,x,y,z,\alpha)|&\leq k_f((1+|x|^\rho+|x'|^\rho)|x'-x|+|y'-y|+|z'-z|).
  \end{align*}
  Moreover, we assume that $f(t,x,y,z,\cdot)$ is continuous for all $(t,x,y,z)\in [0,T]\times\R^n\times\R\times\R^d\to\R$.
  \item\label{ass:on-coeff-psi} The terminal reward $\psi:\R^n\to\R$ satisfies the growth condition
  \begin{align*}
    |\psi(x)|\leq C(1+|x|^\rho)
  \end{align*}
  for all $x\in\R^n$, and the following local Lipschitz criterion
  \begin{align*}
    |\psi(x)-\psi(x')|\leq C(1+|x|^\rho+|x'|^\rho)|x-x'|.
  \end{align*}
  \item\label{ass:on-coeff-ell} The intervention cost $\ell:[0,T]\times \R^n\times U\to \R_+$ is jointly continuous in $(t,x,b)$, bounded from below, \ie
  \begin{align*}
    \ell(t,x,b)\geq\delta >0,
  \end{align*}
  locally Lipschitz in $x$ and locally H\"older continuous in $t$, in particular, we assume that
  \begin{align*}
    |\ell(t,x,b)-\ell(t',x',b)|\leq C(1+|x'|^\rho+|x|^\rho)(|x-x'|+|t'-t|^\varsigma),
  \end{align*}
  for some $\varsigma>0$.
  \item\label{ass:on-coeff-@end} For each $(x,b)\in\R^n\times U$ we have
  \begin{align*}
    \psi(x)>\psi(\Gamma(T,x,b))-\ell(t,x,b).
  \end{align*}
\end{enumerate}
\end{ass}

\begin{rem}\label{rem:@end}
Note in particular that Assumption~\ref{ass:on-coeff}.\ref{ass:on-coeff-@end} implies that the lower and upper value functions defined in the introduction satisfies $\Vlow(T,x)=\Vup(T,x)=\psi(T,x)$ for all $x\in\R^n$.
\end{rem}

Moreover, we make the following assumptions on the coefficients of the controlled forward SDE:

\begin{ass}\label{ass:onSFDE}
For any $t,t'\in [0,T]$, $b\in U$, $\alpha\in A$ and $x,x'\in\R^n$ we have:
\begin{enumerate}[i)]
  \item\label{ass:onSFDE-Gamma} The function $\Gamma:[0,T]\times\R^n\times U\to\R^d$ is jointly continuous and satisfies
  \begin{align*}
    |\Gamma(t,x,b)-\Gamma(t',x',b)|&\leq k_{\Gamma}(|x'-x|+|t'-t|^\varsigma(1+|x|+|x'|))
  \end{align*}
  and the growth condition
  \begin{align}\label{ekv:imp-bound}
    |\Gamma(t,x,b)|\leq K_\Gamma\vee |x|.
  \end{align}
  for some constants $k_\Gamma,K_\Gamma>0$ and $\varsigma>0$.
  \item\label{ass:onSFDE-a-sigma} The coefficients $a:[0,T]\times\R^n\times A\to\R^{n}$ and $\sigma:[0,T]\times\R^n\times A\to\R^{n\times d}$ are jointly continuous and satisfy the growth condition
  \begin{align*}
    |a(t,x,\alpha)|+|\sigma(t,x,\alpha)|&\leq C(1+|x|),
  \end{align*}
  and the Lipschitz continuity
  \begin{align*}
    |a(t,x,\alpha)-a(t,x',\alpha)|+|\sigma(t,x,\alpha)-\sigma(t,x',\alpha)|&\leq C|x'-x|.
  \end{align*}
\end{enumerate}
\end{ass}

\subsection{Viscosity solutions}
We define the upper, $v^*$, and lower, $v_*$ semi-continuous envelope of a function $v$ as
\begin{align*}
v^*(t,x):=\limsup_{(t',x')\to(t,x),\,t'<T}v(t',x')\quad {\rm and}\quad v_*(t,x):=\liminf_{(t',x')\to(t,x),\,t'<T}v(t',x')
\end{align*}
Next we introduce the notion of a viscosity solution using the limiting parabolic superjet $\bar J^+v$ and subjet $\bar J^-v$ of a function $v$ (see pp. 9-10 of \cite{UsersGuide} for a definition):
\begin{defn}\label{def:visc-sol-jets}
Let $v$ be a locally bounded l.s.c. (resp. u.s.c.) function from $[0,T]\times \R^n$ to $\R$. Then,
\begin{enumerate}[a)]
  \item It is referred to as a viscosity supersolution (resp. subsolution) to \eqref{ekv:var-ineq} if:
  \begin{enumerate}[i)]
    \item $v(T,x)\geq \psi(x)$ (resp. $v(T,x)\leq \psi(x)$)
    \item For any $(t,x)\in [0,T)\times\R^d$ and $(p,q,X)\in \bar J^- v(t,x)$ (resp. $J^+ v(t,x)$) we have
    \begin{align*}
      \min\Big\{&v(t,x)-\mcM v(t,x),-p-\inf_{\alpha\in A}H(t,x,v(t,x),q,X,a)\Big\}\geq 0
    \end{align*}
    (resp.
    \begin{align*}
      \min\Big\{&v(t,x)-\mcM v(t,x),-p-\inf_{\alpha\in A}H(t,x,v(t,x),q,X,a)\Big\}\leq 0).
    \end{align*}
  \end{enumerate}
  \item It is referred to as a viscosity solution if it is both a supersolution and a subsolution.
\end{enumerate}
\end{defn}

%\begin{defn}\label{def:visc-sol-jets}
%Let $v$ be a locally bounded l.s.c. (resp. u.s.c.) function from $[0,T]\times \R^n$ to $\R$. Then,
%\begin{enumerate}[a)]
%  \item It is referred to as a viscosity supersolution (resp. subsolution) to \eqref{ekv:var-ineq} if:
%  \begin{enumerate}[i)]
%    \item $v_*(T,x)\geq \psi(x)$ (resp. $v^*(T,x)\leq \psi(x)$)
%    \item For any $(t,x)\in [0,T)\times\R^d$ and $(p,q,X)\in \bar J^- v_*(t,x)$ (resp. $J^+ v^*(t,x)$) we have
%    \begin{align*}
%      \min\Big\{&v_*(t,x)-\mcM v_*(t,x),-p-\inf_{\alpha\in A}H(t,x,v_*(t,x),q,X,a)\Big\}\geq 0
%    \end{align*}
%    (resp.
%    \begin{align*}
%      \min\Big\{&v^*(t,x)-\mcM v^*(t,x),-p-\inf_{\alpha\in A}H(t,x,v^*(t,x),q,X,a)\Big\}\leq 0).
%    \end{align*}
%  \end{enumerate}
%  \item It is referred to as a viscosity solution if it is both a supersolution and a subsolution.
%\end{enumerate}
%\end{defn}

We will sometimes use the following equivalent definition of viscosity supersolutions (resp. subsolutions):
\begin{defn}\label{def:visc-sol-dom}
  A l.s.c.~(resp. u.s.c.) function $v$ is a viscosity supersolution (subsolution) to \eqref{ekv:var-ineq} if $v(T,x)\geq \psi(x)$ (resp. $\leq \psi(x)$) and whenever $\varphi\in C^{3}_{l,b}([0,T]\times\R^d\to\R)$ is such that $\varphi(t,x)=v(t,x)$ and $\varphi-v$ has a local maximum (resp. minimum) at $(t,x)$, then
  \begin{align*}
    \min\big\{&v(t,x)-\mcM v(t,x),-\varphi_t(t,x)-\inf_{\alpha\in A}H(t,x,v(t,x),D\varphi(t,x),D^2\varphi(t,x),a)\big\}\geq 0\:(\leq 0).
  \end{align*}
\end{defn}

\begin{rem}
$C^{3}_{l,b}$ denotes the set of real-valued functions that are continuously differentiable up to third order and whose derivatives of order one to three are bounded
\end{rem}

\subsection{Backward semigroups}
For $(t,x)\in [0,T]\times \R^n$ we let $h\in [0,T-t]$ and assume that $\eta\in L^2(\Omega,\mcF_{t+h},\Prob)$. For all $(u,\alpha)\in\mcU_{t,t+h}\times\mcA_{t,t+h}$ we then define (see \cite{LiPeng09})
\begin{align}
G_{t,t+h}^{t,x;u,\alpha}[\eta]:=\mcY_t,
\end{align}
where $(\mcY,\mcZ)\in\mcS^2\times\mcH^2$ is the unique solution\footnote{From now on we assume that any referred to uniqueness of solutions to a BSDE is uniqueness in $\mcS^2\times\mcH^2$ and therefore refrain from referring to the space.} to
\begin{align*}
\mcY_s&=\eta+\int_s^{t+h}f(r,X^{t,x;u,\alpha}_r,\mcY_r,\mcZ_r)dr-\int_s^{t+h} \mcZ_rdW_r
-\Xi^{t,x;u,\alpha}_{(t+h)+}+\Xi^{t,x;u,\alpha}_s.
\end{align*}
The so defined family of operators $G^{t,x;u,\alpha}$ is referred to as the backward semigroup related to the BSDE.

We note that by the uniqueness of solutions to \eqref{ekv:non-ref-bsde} (see the next section) we have that
\begin{align}\label{ekv:G-is-semigroup}
G_{t,T}^{t,x;u,\alpha}[\psi(X^{t,x;u,\alpha}_T)]=G_{t,t+h}^{t,x;u_{t+h},\alpha}[Y^{t+h,X^{t,x;u_{t+h},\alpha}_{t+h};u^{t+h},\alpha}_{t+h}].
\end{align}
We refer to \eqref{ekv:G-is-semigroup} as the semigroup property of $G$.

%Backward semigroups were introduced by Peng in ...
\section{Forward- Backward SDEs with impulses\label{sec:FBSDEs}}
In this section we consider the non-standard BSDE in \eqref{ekv:non-ref-bsde}. Impulsively controlled BSDEs in the non-Markovian framework were treated in \cite{rbsde_impulse}, while BSDEs related to switching problems have been treated in \cite{HuTang,HamZhang,Morlais13}.

Considering first the forward SDE, we get by repeated use of standard results for SDEs (see \eg Chapter 5 in \cite{Protter}) that \eqref{ekv:forward-sde1}-\eqref{ekv:forward-sde2} admits a unique solution $X^{t,x;u,\alpha}$ for any $(u,\alpha)\in\mcU\times\mcA$ since $N<\infty$, $\Prob$-a.s. Now, any solution of \eqref{ekv:non-ref-bsde} can be written $Y^{t,x;u,\alpha}_s=\tilde Y^{t,x;u,\alpha}_s+\Xi^{t,x;u,\alpha}_s$, where $(\tilde Y^{t,x;u,\alpha},\tilde Z^{t,x;u,\alpha})\in \mcS^2_c\times\mcH$ solves the standard BSDE
\begin{align}
\tilde Y^{t,x;u,\alpha}_s&=\psi(X^{t,x;u,\alpha}_T)-\Xi^{t,x;u,\alpha}_{T+}+\int_s^Tf(r,X^{t,x;u,\alpha}_r,\tilde Y^{t,x;u,\alpha}_r+\Xi^{t,x;u,\alpha}_r,\tilde Z^{t,x;u,\alpha}_r)dr-\int_s^T \tilde Z^{t,x;u,\alpha}_rdW_r. \label{ekv:non-ref-bsde-std}
\end{align}
By standard results we find that \eqref{ekv:non-ref-bsde-std} admits a unique solution whenever $\Xi^{t,x;u,\alpha}_{T+}\in L^2(\Omega,\Prob)$ and $f(\cdot,X^{t,x;u,\alpha}_\cdot,0,0)\in\mcH^2$. By a moment estimate given in the next section we are able to conclude that \eqref{ekv:non-ref-bsde} admits a unique solution whenever $(u,\alpha)\in\mcU\times\mcA$.

%In this section we give estimates for solutions of \eqref{ekv:non-ref-bsde}. To obtain such estimates we first need to derive moment and stability estimates for solutions to our controlled SDE \eqref{ekv:forward-sde1}-\eqref{ekv:forward-sde2}.

\subsection{Estimates for the controlled diffusion process}
%As mentioned above, the following moment estimate guarantees existence of a unique solution to the BSDE \eqref{ekv:non-ref-bsde}.
\begin{prop}\label{prop:SDEmoment}
%Under Assumption~\ref{ass:on-coeff},
For each $p\geq 1$, there is a $C>0$ such that
\begin{align}\label{ekv:SDEmoment}
\E\Big[\sup_{s\in[\zeta,T]}|X^{t,x;u,\alpha}_s|^{p}\Big|\mcF_\zeta\Big]\leq C(1+|X^{t,x;u,\alpha}_\zeta|^{p}),
\end{align}
$\Prob$-a.s.~for all $(t,\zeta,x,u,\alpha)\in [0,T]^2\times\R^n\times\mcU\times\mcA$.
\end{prop}

\noindent\emph{Proof.} We use the shorthand $X^j:=X^{t,x;[u]_j,\alpha}$. By Assumption~\ref{ass:onSFDE}.(\ref{ass:onSFDE-Gamma}) we get for $s\in [\tau_{j},T]$, using integration by parts, that
\begin{align*}
|X^{j}_s|^2&= |X^{j}_{\tau_{j}}|^2+2\int_{\tau_{j}+}^s X^{j}_{r}dX^{j}_r+\int_{\tau_{j}+}^s d[X^{j},X^{j}]_r
\\
&\leq K^2_\Gamma\vee |X^{{j-1}}_{\tau_{j}}|^2+2\int_{\tau_{j}+}^s X^{j}_{r} dX^{j}_r+\int_{\tau_{j}+}^s d[X^{j},X^{j}]_r.
\end{align*}
We note that if $|X^{j}_s|> K_\Gamma$ and $|X^{j}_r|\leq K_\Gamma$ for some $r\in[\zeta,s)$ then there is a largest time $\theta<s$ such that $|X^{j}_{\theta}|\leq K_\Gamma$. This means that during the interval $(\theta,s]$ interventions will not increase the magnitude $|X^{j}|$. By induction we find that
\begin{align}\label{ekv:X2-bound}
|X^{j}_s|^2&\leq |X^{j}_\zeta|^2\vee K_\Gamma^2+\sum_{i=0}^{j} \Big\{2\int_{\theta\vee(\tilde\tau_{i}+)}^{s\wedge\tilde\tau_{i+1}} X^{i}_{r}dX^{i}_r+\int_{\theta\vee(\tilde\tau_{i}+)}^{s\wedge\tilde\tau_{i+1}} d[X^{i},X^{i}]_r\Big\}
\end{align}
for all $s\in[t,T]$, where $\theta:=\sup\{r\geq 0 : |X^u_r|\leq K_\Gamma\}\vee \zeta$, $\tilde\tau_0+=0$, $\tilde\tau_i=\tau_i$ for $i=1,\ldots,j$ and $\tilde\tau_{j+1}=\infty$.

Now, since $X^{i}$ and $X^{j}$ coincide on $[0,\tau_{i+1\wedge j+1})$ we have
\begin{align*}
\sum_{i=0}^{j}\int_{\theta\vee\tilde\tau_{i}+}^{s\wedge\tilde\tau_{i+1}} X^{i}_{r} dX^{i}_r
&=\int_{\theta}^s X^{j}_{r}a(r,X^{j}_r,\alpha_r)dr+\int_{\theta}^{s}X^{j}_{r}\sigma(r,X^{j}_r,\alpha_r)dW_r,
\end{align*}
and
\begin{align*}
\sum_{i=0}^{j} \int_{\theta\vee\tilde\tau_{i}+}^{s\wedge\tilde\tau_{i+1}} d[X^{i},X^{i}]_r&=\int_{\theta}^{s} \sigma^2(r,X^{j}_r,\alpha_r)dr.
\end{align*}
Inserted in \eqref{ekv:X2-bound} this gives
\begin{align*}
|X^{j}_s|^2&\leq |X^{j}_\zeta|^2\vee K_\Gamma^2+\int_{\theta}^s (2X^{j}_{s}a(r,X^{j}_r,\alpha_r)+\sigma^2(r,X^{j}_r,\alpha_r))dr+2\int_{\theta}^{s}X^{j}_{r}\sigma(r,X^{j}_r,\alpha_r)dW_r
\\
&\leq |X^{j}_\zeta|^2+C\Big(1+\int_{\zeta}^{s}|X^{j}_{r}|^2dr+\sup_{v\in[\zeta,s]}\Big|\int_{\zeta}^{v}X^{j}_r\sigma(r,X^{j}_r)dW_r\Big|\Big).
\end{align*}
The Burkholder-Davis-Gundy inequality now gives that for $p\geq 2$,
\begin{align*}
\E\Big[\sup_{r\in[\zeta,s]}|X^{j}_r|^{p}\Big|\mcF_\zeta\Big]\leq |X^{j}_{\zeta}|^p +C\big(1+\E\Big[\int_{\zeta}^{s}|X^{i}_{r}|^{p}dr+\big(\int_{\zeta}^{s}|X^{j}_r|^4 dr\big)^{p/4}\Big]\big)
\end{align*}
and Gr\"onwall's lemma gives that for $p\geq 4$,
\begin{align}\label{ekv:moment_steg1}
\E\Big[\sup_{s\in[\zeta,T]}|X^{j}_s|^{p}\big|\mcF_\zeta\Big]&\leq C(1+ |X^{j}_\zeta|^{p}),
\end{align}
$\Prob$-a.s., where the constant $C=C(T,p)$ does not depend on $u$, $\alpha$ or $j$ and \eqref{ekv:SDEmoment} follows by letting $j\to\infty$ on both sides and using Fatou's lemma. The result for general $p\geq 1$ follows by Jensen's inequality.\qed\\

As mentioned above, inequality \eqref{ekv:SDEmoment} guarantees existence of a unique solution to the BSDE \eqref{ekv:non-ref-bsde}. We will also need the following stability property.
\begin{prop}\label{prop:SFDEflow}
For each $k\geq 0$ and $p\geq 1$, there is a $C\geq 0$ such that
\begin{align*}
  \E\Big[\sup_{s\in[t',T]}|X^{t,x;u,\alpha}_{s} -X^{t',x';u,\alpha}_{s}|^{p}\Big|\mcF_t\Big]\leq C(|x-x'|^p+(1+|x|^p)|t'-t|^{p(\varsigma\wedge 1/2)}),
\end{align*}
$\Prob$-a.s.~for all $(t,t',x,x')\in [0,T]^2\times \R^{2n}$, with $t'\geq t$, and all $(u,\alpha)\in\mcU^k\times\mcA$.
\end{prop}

\noindent\emph{Proof.} To simplify notation we let $X^j:=X^{t,x;[u]_j,\alpha}$ and $X^{'j}:=X^{t',x';[u]_j,\alpha}$ for $j=0,\ldots,k$. Moreover, we let $\delta X^j:=X^j- X^{'j}$ and set $\delta X:=\delta  X^{k}$. Define $\kappa:=\max\{j\geq 0:\tau_j\leq t'\}\vee 0$, then if $\kappa=0$ we have $|\delta  X_{t'}|=|\delta  X^0_{t'}|$, where for any value of $\kappa$,
\begin{align*}
|\delta  X^0_{t'}|=|X^{t,x;u,\alpha}_{t'}-x'|.
\end{align*}
When $\kappa>0$ we get for $j=1,\ldots,\kappa$,
\begin{align*}
|\delta  X^{j}_{t'}|&\leq k_\Gamma(|\delta X^{{j-1}}_{t'}|
+|X^{{j-1}}_{t'}-X^{{j-1}}_{\tau_j}|+|t'-t|^\varsigma(1+|X^{{j-1}}_{\tau_j}|+|X^{'{j-1}}_{t'}|))
\end{align*}
By induction we find that
\begin{align*}
|\delta  X^{\kappa}_{t'}|&\leq \sum_{j=1}^{\kappa} k_\Gamma^{\kappa+1-j}(|X^{{j-1}}_{t'}-X^{{j-1}}_{\tau_{j}}| + |t'-t|^\varsigma(1+\sup_{s\in [t,t']}|X^{{j-1}}_{s}|+|X^{'{j-1}}_{t'}|)).
\end{align*}
Now, since
\begin{align*}
|X^{{j-1}}_{t'}-X^{{j-1}}_{\tau_j}|\leq \int_{\tau_j}^{t'}|a(r,X^{j-1}_r,\alpha_r)|dr+\Big| \int_{\tau_j}^{t'}\sigma(r,X^{j-1}_r,\alpha_r)dW_s\Big|,
\end{align*}
Proposition~\ref{prop:BSDEmoment} gives that
\begin{align*}
\E\big[|X^{j-1}_{t'}-X^{j-1}_{\tau_j}|^p\big|\mcF_t\big] \leq C(1+|x|^p)|t'-t|^{p/2}.
\end{align*}
Similarly,
\begin{align*}
 \E\big[|X^{t,x;u,\alpha}_{t'}-x'|^p\big|\mcF_t\big] \leq C(|x-x'|^p+(1+|x|^p)|t'-t|^{p/2})
\end{align*}
and we find that
\begin{align*}
\E\big[|\delta  X_{t'}|^p\big|\mcF_t\big] \leq C(|x-x'|^p+(1+|x|^p)|t'-t|^{p(\varsigma\wedge 1/2)}).
\end{align*}

Moreover, we note that for $j\geq\kappa$ and $s\geq \tau_j$ (with $\tau_0:=t'$),
\begin{align*}
|\delta  X^{j}_{s}|&\leq (1\vee k_\Gamma)|\delta  X^j_{\tau_j}|+\int_{\tau_j}^{s}|a(r,X^{j}_r,\alpha_r)-a(r, X^{'j}_r,\alpha_r)|dr
\\
&\quad+\Big|\int_{t'}^{s}(\sigma(r,X^{j}_r,\alpha_r)-\sigma(r,X^{'j}_r,\alpha_r)dW_r\Big|
\end{align*}
and the Burkholder-Davis-Gundy inequality gives for $p\geq 2$ we have
\begin{align*}
\E\Big[\sup_{r\in[\tau_j,s]}|\delta X^j_{r}|^{p}\Big] &\leq C\E\Big[|\delta X^j_{\tau_j}|^{p}+ \Big(\int_{\tau_j}^s|a(r,X^{j}_r,\alpha_r)- a(r, X^{'j}_r,\alpha_r)|dr\Big)^{p}
\\
& + \Big(\int_{\tau_j}^{s}|\sigma(r, X^{'j}_r,\alpha_r)-\sigma(r,X^j_r,\alpha_r)|^2dr\Big)^{p/2}\Big]
\\
&\leq C\E\Big[|\delta X^j_{\tau_j}|^{p}+ \big(\int_{\tau_j}^s|\delta X^j_{r}|^{2}dr\big)^{p/2}\Big].
\end{align*}
The Lipschitz conditions on the coefficients combined with Gr\"onwall's lemma then implies that
\begin{align*}
\E\Big[\sup_{r\in[\tau_j,T]}|\delta X^j_{r}|^{p}\Big] &\leq C\E\Big[|\delta X^j_{\tau_j}|^{p}\Big].
\end{align*}
Now, since $|\delta X^l_{\tau_l}|\leq k_\Gamma|\delta X^{l-1}_{\tau_l}|$ for $l=\kappa+1,\ldots,N$ the result follows by induction.\qed\\

\subsection{Estimates for the BSDE}
For $(t,x)\in[0,T]\times\R^n$ and $(u,\alpha)\in\mcU\times\mcA$ we let $(\check Y^{t,x;u,\alpha},\check Z^{t,x;u,\alpha})$ be the unique solution to the following standard BSDE%\in \mcS^2\times\mcH^2
\begin{align}\label{ekv:bsde-trad}
  \check Y_s^{t,x;u,\alpha}=\psi(X^{t,x;u,\alpha}_T)+\int_s^Tf(r,X^{t,x;u,\alpha}_r,\check Y^{t,x;u,\alpha}_r,\check Z^{t,x;u,\alpha}_r)dr-\int_s^T \check Z^{t,x;u,\alpha}_rdW_r.
\end{align}
Combining classical results (see \eg \cite{ElKaroui2}) with Proposition~\ref{prop:SDEmoment}, we have
\begin{align}\nonumber
  &\E\Big[\sup_{s\in [t,T]}|\check Y_s^{t,x;u,\alpha}|^2+\int_t^T|\check Z^{t,x;u,\alpha}_s|^2ds\Big|\mcF_t\Big]
  \\
  &\leq C\E\Big[|\psi(X^{t,x;u,\alpha}_T)|^2+\int_t^T|f(r,X^{t,x;u,\alpha}_r,0,0,\alpha_r)|^2dr\Big|\mcF_t\Big]\leq C(1+|x|^{2\rho}),\label{ekv:bsde-trad-moment}
\end{align}
$\Prob$-a.s.~for all $(u,\alpha)\in\mcU\times\mcA$.

We have the following straightforward generalization of the standard comparison principle:
\begin{lem} (Comparison principle)
If $\hat f$ satisfies Assumption~\ref{ass:on-coeff}, and $\hat G^{t,x;u,\alpha}$ is defined as $G^{t,x;u,\alpha}$ but with driver $\hat f$ instead of $f$, then if $f(t,x,y,z,\alpha)\leq \hat f(t,x,y,z,\alpha)$ for all $(t,x,y,z,\alpha)\in [0,T]\times\R^d\times\R\times\R^d\times U$, we have $G_{s,r}^{t,x;u,\alpha}[\eta]\leq \hat G_{s,r}^{t,x;u,\alpha}[\hat\eta]$, $\Prob$-a.s.~for each $t\leq s\leq r\leq T$ whenever $\eta,\hat\eta\in L^2(\Omega,\mcF_s,\Prob)$ are such that $\eta\leq \hat\eta$, $\Prob$-a.s.
\end{lem}

\noindent\emph{Proof.} This follows immediately from the standard comparison principle (see Theorem 2.2 in \cite{ElKaroui2}).\qed\\

%\begin{enumerate}[a)]
%\item
%\item If, in addition $\Prob[\hat\psi(X^{t,x;u,\alpha}_T)>\psi(X^{t,x;u,\alpha}_T)]>0$, then $\Prob[\hat Y_s^{t,x;u,\alpha}> Y_s^{t,x;u,\alpha}]>0$.
%\end{enumerate}

Using the comparison principle we easily deduce the following moment estimates:
\begin{prop}\label{prop:BSDEmoment}
%Under Assumption~\ref{ass:on-coeff},
We have,
\begin{align}\label{ekv:Ybnd}
  \esssup_{\alpha\in\mcA}|\esssup_{u\in\mcU}Y^{t,x;u,\alpha}_t|\leq C(1+|x|^\rho),\qquad \Prob-{\rm a.s.}
\end{align}
and for each $k\geq 0$, there is a $C>0$ such that
\begin{align}\label{ekv:BSDEmoment}
\E\Big[\sup_{s\in[t,T]}|Y^{t,x;u,\alpha}_s|^{2}+\int_t^T| Z^{t,x;u,\alpha}_s|^2ds\Big|\mcF_t\Big]\leq C(1+|x|^{2\rho}),
\end{align}
$\Prob$-a.s.~for all $(t,x,u,\alpha)\in [0,T]\times\R^n\times\mcU^k\times\mcA$.
\end{prop}

\noindent\emph{Proof.} The first statement follows by repeated application of the comparison principle which gives that $\check Y^{t,x;\emptyset,\alpha}_t\leq\esssup_{u\in\mcU}Y^{t,x;u,\alpha}_t \leq \esssup_{u\in\mcU}\check Y^{t,x;u,\alpha}_t$ and using~\eqref{ekv:bsde-trad-moment}.

The second statement follows by noting that for fixed $k\geq 0$, there is a $C>0$ such that
\begin{align*}
\E[|\Xi^{t,x;u,\alpha}_{T+}|^2]&\leq C(1+\E[\sup_{s\in[t,T]}|X^{t,x;u,\alpha}_{s}|^{2\rho}])\leq C(1+|x|^{2\rho})
\end{align*}
for all $(u,\alpha)\in \mcU^k\times\mcA$.\qed\\

\begin{prop}\label{prop:BSDEflow}
For each $k\geq 0$, there is a $C>0$ such that
\begin{align}\label{ekv:BSDEmoment}
|\E\big[Y^{t',x';u,\alpha}_{t'}-Y^{t,x;u,\alpha}_{t}\big]|\leq C(1+|x|^{\rho+1}+|x'|^{\rho+1})(|x'-x|+|t'-t|^{\varsigma\wedge 1/2}),
\end{align}
$\Prob$-a.s.~for all $(t,x),(t',x')\in[0,T]\times\R^n$ with $t\leq t'$ and all $u\in\mcU^k$ and $\alpha\in\mcA$.
\end{prop}

\noindent\emph{Proof.}  To simplify notation, we let $X:=X^{t,x;u,\alpha}$ and $X':=X^{t',x';u,\alpha}$ and set $(Y,Z):=(Y^{t,x;u,\alpha},Z^{t,x;u,\alpha})$ and $(Y',Z'):=(Y^{t',x';u,\alpha},Z^{t',x';u,\alpha})$. By defining $\delta Y:=Y-Y'_{\cdot\vee t'}$ and $\delta Z:=Z-\ett_{[\cdot\geq t']}Z'$ we have for $s\in[t,T]$ that
\begin{align*}
\delta Y_{s}&=\psi(X_T)-\psi(X'_T)+\int_s^{T}(f(r,X_r,Y_r,Z_r,\alpha_r)-f(r,X'_r,Y'_r,Z'_r,\alpha_r))dr
\\
&\quad-\int_s^{T} \delta Z_rdW_r-\sum_{j=1}^N(\ett_{[\tau_j \geq s]}\ell(\tau_j,X^{j-1}_{\tau_j},\beta_j)-\ett_{[\tau_j\vee t' \geq s]}\ell(\tau_j\vee t',X^{'j-1}_{\tau_j\vee t'},\beta_j)),
\end{align*}
with $X^j:=X^{t,x;[u]_j,\alpha}$ and $X^{'j}:=X^{t',x';[u]_j,\alpha}$. We now introduce the processes $(\zeta_1(s))_{s\in [t,T]}$ and $(\zeta_2(s))_{s\in [t,T]}$ defined as\footnote{Throughout, we use the convention that $\frac{0}{0}0=0$}
\begin{align*}
 \zeta_1(s):=\frac{f(s,X_s,Y_s,Z_s,\alpha_s)-f(s,X_s,\ett_{[s\geq t']}Y'_s,Z_s,\alpha_s)}{Y_s-\ett_{[s\geq t']}Y'_s}\ett_{[Y_s\neq \ett_{[s\geq t']}Y'_s]}
\end{align*}
and
\begin{align*}
 \zeta_2(s):=\frac{f(s,X_s,\ett_{[s\geq t']}Y'_s,Z_s,\alpha_s)-f(s,X_s,\ett_{[s\geq t']}Y'_s,\ett_{[s\geq t']}Z'_s,\alpha_s)}{|Z_s-\ett_{[s\geq t']}Z'_s|^2}(Z_s-\ett_{[s\geq t']}Z'_s)^\top.
\end{align*}
We then have by the Lipschitz continuity of $f$ that $|\zeta_1(s)|\vee|\zeta_2(s)|\leq k_f$. Using Ito's formula we find that
\begin{align*}
\delta Y_s&=R_{s,T}(\psi(X_T)-\psi(X'_T))+\int_s^{T}R_{s,r}(f(r,X_r,\ett_{[r\geq t']}Y'_r,\ett_{[r\geq t']}Z'_r,\alpha_r)-\ett_{[r\geq t']}f(r,X'_r,Y'_r,Z'_r,\alpha_r))dr
\\
&\quad-\int_s^{T} R_{s,r}\delta Z_rdW_r-\sum_{j=1}^N(\ett_{[\tau_j \geq s]}R_{s,\tau_j}\ell(\tau_j,X^{j-1}_{\tau_j},\beta_j)-\ett_{[\tau_j\vee t' \geq s]}R_{s,\tau_j\vee t'}\ell(\tau_j\vee t',X^{'j-1}_{\tau_j\vee t'},\beta_j))
\end{align*}
with $R_{s,r}:=e^{\int_s^{r}(\zeta_1(v)-\frac{1}{2}|\zeta_2(v)|^2)dv+\frac{1}{2}\int_s^{r}\zeta^u_2(v)dW_v}$. Taking expectations on both sides yields
\begin{align*}
  |\E\big[\delta Y_t\big]|&\leq C\E\Big[R_{t,T}(1+|X_T|^\rho+|X'_T|^\rho)|X'_T-X_T|+\int_t^{t'}R_{t,r}(1+|X_r|^\rho)dr
  \\
  &\quad+\int_{t'}^{T}R_{t,r}(1+|X_r|^\rho+|X'_r|^\rho)|X'_r-X_r|dr
\\
&\quad+\sum_{j=1}^N R_{t,\tau_j}|\ell(\tau_j,X^{j-1}_{\tau_j},\beta_j)-R_{\tau_j,\tau_j\vee t'}\ell(\tau_j\vee t',X^{'j-1}_{\tau_j\vee t'},\beta_j)|\Big].
\end{align*}
Now,
\begin{align*}
  &\E\Big[R_{t,T}(1+|X_T|^\rho+|X'_T|^\rho)|X'_T-X_T|+\int_t^{t'}R_{t,r}(1+|X_r|^\rho)dr+\int_{t'}^{T}R_{t,r}(1+|X_r|^\rho+|X'_r|^\rho)|X'_r-X_r|dr\Big]
  \\
  &\leq C\E\Big[\sup_{s\in [t,T]}|R_{t,s}|^2\Big]^{1/2}\E\Big[(t'-t)\int_t^{t'}(1+|X_r|^{2\rho})dr+\sup_{r\in [t',T]}(1+|X_r|^{2\rho}+|X'_r|^{2\rho})|X'_r-X_r|^2\Big]^{1/2}
  \\
  &\leq C(|t'-t|+\E\Big[\sup_{r\in [t',T]}(1+|X_r|^{4\rho}+|X'_r|^{4\rho})\Big]^{1/4}\E\Big[\sup_{r\in [t',T]}|X'_r-X_r|^4\Big]^{1/4})
  \\
  &\leq C(1+|x|^\rho+|x'|^\rho)(|x-x'|+(1+|x|)|t'-t|^{\varsigma\wedge 1/2})
\end{align*}
where we have used Proposition~\ref{prop:SFDEflow} to reach the last inequality. Moreover,
\begin{align*}
  &\E\Big[\sum_{j=1}^N R_{t,\tau_j}|\ell(\tau_j,X^{j-1}_{\tau_j},\beta_j)-R_{\tau_j,\tau_j\vee t'}\ell(\tau_j\vee t',X^{'j-1}_{\tau_j\vee t'},\beta_j)|\Big]
  \\
  &\leq  \E\Big[\sum_{j=1}^N R_{t,\tau_j}\big((1+R_{\tau_j,\tau_j\vee t'})|\ell(\tau_j,X^{j-1}_{\tau_j},\beta_j)-\ell(\tau_j\vee t',X^{'j-1}_{\tau_j\vee t'},\beta_j)|
  \\
  &\quad+|1-R_{\tau_j,\tau_j\vee t'}|(\ell(\tau_j,X^{j-1}_{\tau_j},\beta_j)+\ell(\tau_j\vee t',X^{'j-1}_{\tau_j\vee t'},\beta_j))\big)\Big]
  \\
  &\leq Ck\E\Big[\sup_{s\in [t,T]}|R_{t,s}|^2\Big]^{1/2}\Big(\E\Big[\sup_{r\in [t,T]}(1+|X_r|^{2\rho}+|X'_r|^{2\rho})(|X'_{r\vee t'}-X_r|^2+|t'-t|^{2\varsigma})\Big]^{1/2}
  \\
  &\quad + \E\Big[\sup_{r\in [t,t']}|1-R_{t,r}|^2(1+|X_r|^{2\rho})\Big]^{1/2}
  \\
  &\leq Ck(1+|x|^\rho+|x'|^\rho)(|x'-x|+(1+|x|)|t'-t|^{ \varsigma\wedge1/2}).
\end{align*}
Combining the above inequalities, the assertion follows.\qed\\

The above proof immediately gives the following stability result:
\begin{cor} (Stability)
If $\hat f$ satisfies Assumption~\ref{ass:on-coeff}, and $\hat G^{t,x;u,\alpha}$ is defined as $G^{t,x;u,\alpha}$ with driver $\hat f$ instead of $f$, then there is a $C>0$ such that
\begin{align*}
|\hat G_{t,s}^{t,x;u,\alpha}[\hat\eta]-G_{t,s}^{t,x;u,\alpha}[\eta]|\leq C\E\Big[|\hat\eta-\eta|^2+\int_t^{s}|\hat f(r,X_r,\mcY_r,\mcZ_r,\alpha_r)-f(r,X_r,\mcY_r,\mcZ_r,\alpha_r)|^2dr \Big|\mcF_t\Big]^{1/2},
\end{align*}
$\Prob$-a.s.~for all $s\in [t,T]$ and $\eta,\hat\eta\in L^2(\Omega,\mcF_{s},\Prob)$.
\end{cor}

\section{Dynamic programming principles\label{sec:DPP}}
In this section we show that $\Vlow$ and $\Vup$ are jointly continuous (deterministic) functions that satisfy the dynamic programming relations
\begin{align}
\Vlow(t,x)=\essinf_{\alpha^S\in\mcA^S_{t,t+h}}\esssup_{u\in\mcU_{t,t+h}}G_{t,t+h}^{t,x;u,\alpha^S(u)}[ \Vlow(t+h,X^{t,x;u,\alpha^S(u)}_{t+h})]\label{ekv:dynp-W}
\end{align}
and
\begin{align}
\Vup(t,x)=\esssup_{u^S\in\mcU^S_{t,t+h}}\essinf_{\alpha\in\mcA_{t,t+h}}G_{t,t+h}^{t,x;u^S(\alpha),\alpha}[ \Vup(t+h,X^{t,x;u^S(\alpha),\alpha}_{t+h})],\label{ekv:dynp-U}
\end{align}
for $t\in[0,T]$ and $h\in[0,T-t]$.

\begin{prop}
For every $(t,x)\in[0,T]\times \R^n$ we have $\Vlow(t,x)=\E[\Vlow(t,x)]$ and $\Vup(t,x)=\E[\Vup(t,x)]$, $\Prob$-a.s.
\end{prop}

\noindent\emph{Proof.} This follows by repeating the steps in the proof of Proposition 4.1 in \cite{Buckdahn08}.\qed\\

We can thus pick the deterministic versions to represent $\Vlow$ and $\Vup$. As mentioned in the introduction, the main technical difficult that we encounter appears when trying to show continuity of the upper and lower value functions in the time variable. The reason for this is that the constant $C$ in Proposition~\ref{prop:BSDEflow} depends on $k$ and tends to infinity as $k$ tends to infinity. We resolve this issue by first considering the upper and lower value functions under an imposed restriction on the number of interventions in the impulse control. Relying on a uniform convergence result will then give us continuity of $\Vlow$ and $\Vup$.

\subsection{A DPP with limited number of impulses}
We introduce the truncated value functions
\begin{align*}
  \Vlow^k(t,x):=\essinf_{\alpha^S\in\mcA^S_t}\esssup_{u\in\mcU^k_t}J(t,x;u,\alpha^S(u))
\end{align*}
and
\begin{align*}
  \Vup^k(t,x):=\esssup_{u^S\in\mcU^{S,k}_t}\essinf_{\alpha\in\mcA_t}J(t,x;u^S(\alpha),\alpha)
\end{align*}
for $k\geq 0$. Similarly to $\Vlow$ and $\Vup$ we have:

\begin{lem}
For every $(t,x)\in[0,T]\times \R^n$ and $k\geq 0$ we have $\Vlow^k(t,x)=\E[\Vlow^k(t,x)]$ and $\Vup^k(t,x)=\E[\Vup^k(t,x)]$, $\Prob$-a.s.
\end{lem}

Combined with the estimates of the previous section this gives the following estimates:
\begin{prop}\label{prop:W-k-cont}
For each $k\geq 0$, there is a $C>0$ such that
\begin{align}
  |\Vlow^k(t,x)-\Vlow^k(t',x')|+|\Vup^k(t,x)-\Vup^k(t,x')|\leq C(1+|x|^{\rho+1}+|x|^{\rho+1})(|x'-x|+|t-t'|^{\varsigma\wedge 1/2}),
\end{align}
for all $(t,x),(t',x')\in [0,T]\times\R^n$. Moreover, there is a $C>0$ such that
\begin{align*}
  |\Vlow^k(t,x)|+|\Vup^k(t,x)|\leq C(1+|x|^\rho)
\end{align*}
for all $k\geq 0$ and $(t,x)\in[0,T]\times \R^n$.
\end{prop}

\noindent\emph{Proof.} Since
\begin{align*}
  \Vlow^k(t,x)=\essinf_{\alpha^S\in\mcA^S}\esssup_{u\in\mcU^k}Y^{t,x;u,\alpha^S(u)},
\end{align*}
we have
\begin{align*}
  \Vlow^k(t,x)-\Vlow^k(t',x')&= \essinf_{\alpha^S\in\mcA^S}\esssup_{u\in\mcU^k}Y^{t,x;u,\alpha^S(u)}_t-\essinf_{\alpha^S\in\mcA^S}\esssup_{u\in\mcU^k}Y^{t',x';u,\alpha^S(u)}_{t'}
  \\
  &\leq\esssup_{\alpha^S\in\mcA^S}\{\esssup_{u\in\mcU^k}Y^{t,x;u,\alpha^S(u)}_t-\esssup_{u\in\mcU^k}Y^{t',x';u,\alpha^S(u)}_{t'}\}
  \\
  &\leq \esssup_{\alpha\in\mcA}\esssup_{u\in\mcU^k}\{Y^{t,x;u,\alpha}_t-Y^{t',x';u,\alpha}_{t'}\}
  \\
  &\leq Y^{t,x;u_\epsilon,\alpha_\epsilon}_t-Y^{t',x';u_\epsilon,\alpha_\epsilon}_{t'}+\epsilon
\end{align*}
for each $\epsilon>0$ and some $(u_\epsilon,\alpha_\epsilon)\in \mcU^k\times\mcA$. We also see that the same relation holds for $\Vup^k$. Taking expectation on both sides and using that $\Vlow^k$ and $\Vup^k$ are deterministic, the first inequality follows by Proposition~\ref{prop:BSDEflow} since $\epsilon>0$ was arbitrary.

The second inequality is an immediate consequence of Proposition~\ref{prop:BSDEmoment}.\qed\\

Turning now to the dynamic programming principles, that will be obtained by applying arguments similar to those in Section 4 of \cite{Buckdahn08}, we have:
\begin{prop}\label{prop:dynp-trunk}
For each $k\geq 0$ and any $t\in[0,T]$, $h\in[0,T-t]$ and $x\in \R^n$ we have
\begin{align}\label{ekv:dynp-W-trunk}
  \Vlow^k(t,x)=\essinf_{\alpha^S\in\mcA^S_{t,t+h}}\esssup_{u\in\mcU^k_{t,t+h}}G_{t,t+h}^{t,x;u,\alpha^S(u)}[\Vlow^{k-N}(t+h,X^{t,x;u,\alpha^S(u)}_{t+h})]
\end{align}
and
\begin{align}\label{ekv:dynp-U-trunk}
  \Vup^k(t,x)=\esssup_{u^S\in\mcU^{S,k}_{t,t+h}}\essinf_{\alpha\in\mcA_{t,t+h}}G_{t,t+h}^{t,x;u^S(\alpha),\alpha}[\Vup^{k-N}(t + h,X^{t,x;u^S(\alpha),\alpha}_{t+h})].
\end{align}
\end{prop}

\begin{rem}\label{rem:@tph}
At first glance the DPP for $\Vup$ may seem counter-intuitive as, on the right-hand side, $\alpha$ could take two different values at time $t+h$ (one under $G$ and the other in $\Vup^{k-N}(t + h,\cdot)$) and thus trigger two different reactions from the impulse controller at time $t+h$. However, by the definition of a non-anticipative strategy, $u^S(\alpha)=u^S(\tilde\alpha)$ whenever $\alpha=\tilde\alpha$, $d\Prob\times d\lambda$-a.s.~and an arbitrary choice of $\alpha_{t+h}$ will not influence the overall value.
\end{rem}

\noindent\emph{Proof.} The proof (which is only given for the lower value function $\Vlow^k$ as the arguments for $\Vup^k$ are identical) will be carried out over a sequence of lemmata where
\begin{align*}
  \Vlowh^k(t,x):=\essinf_{\alpha^S\in\mcA^S_{t,t+h}}\esssup_{u\in\mcU^k_{t,t+h}}G_{t,t+h}^{t,x;u,\alpha^S(u)}[ \Vlow^{k-N}(t+h,X^{t,x;u,\alpha^S(u)}_{t+h})].
\end{align*}

\begin{lem}
$\Vlowh^k$ can be chosen to be deterministic.
\end{lem}

\noindent\emph{Proof.} Again, this follows by repeating the steps in the proof of Proposition 4.1 in \cite{Buckdahn08}.\qed\\

\begin{lem}\label{lem:Wh-lessthan-W}
$\Vlowh^k(t,x)\leq \Vlow^k(t,x)$.
\end{lem}

\noindent\emph{Proof.} We begin by picking an arbitrary $\alpha^S\in\mcA^S_t$ and note that we can define the restriction, $\alpha_1^S$, of $\alpha^S$ to $\mcA^S_{t,t+h}$ as
\begin{align*}
  \alpha^S_1(u_1):=\alpha^S(u_1)\big|_{[t,t+h]},\qquad \forall u_1\in\mcU_{t,t+h}.
\end{align*}
%Moreover, we can define the restriction of $\alpha^S$ to $\mcA_{t+h}$ as $\alpha^S_2(u)=\alpha(u)|_{[t+h,T]}$.
We fix $\epsilon>0$ and have by a pasting property\footnote{We can paste together two controls $u_1,u_2\in\mcU^k_s$ on sets $B_1\in\mcF_s$ and $B_2=B_1^c$ by setting $u=\ett_{B_1}u_1+\ett_{B_2}u_2\in\mcU^k_s$ and get by uniqueness of solutions to our BSDE that $G_{s,r}^{t,x;u,\alpha}[\eta]=\ett_{B_1}G_{s,r}^{t,x;u_1,\alpha}[\eta]+\ett_{B_2}G_{s,r}^{t,x;u_2,\alpha}[\eta]$.} %(see the proof of Lemma 4.4 in \cite{Buckdahn08})
that there is a $u^\epsilon_{1}=(\tau^{1,\epsilon}_j,\beta^{1,\epsilon}_j)_{1\leq j\leq N^\epsilon_1}\in \mcU^k_{t,t+h}$ such that
\begin{align*}
\Vlowh^k(t,x)&\leq\esssup_{u\in\mcU^k_{t,t+h}}G_{t,t+h}^{t,x;u,\alpha_1(u)}[\Vlow^{k-N}(t+h,X^{t,x;u_1,\alpha_1(u_1)}_{t+h})]
\\
&\leq G_{t,t+h}^{t,x;u^\epsilon_1,\alpha_1(u^\epsilon_1)}[\Vlow^{k-N^\epsilon_1}(t+h,X^{t,x;u^\epsilon_1,\alpha_1(u^\epsilon_1)}_{t+h})]+\epsilon.
\end{align*}
Now, given $u^\epsilon_{1}$ we can define the restriction, $\alpha^S_2$, of $\alpha^S$ to $\mcA_{t+h}$ as
\begin{align*}
  \alpha^S_2(u_2):=\alpha^S(u_1^\epsilon\oplus u_2)\big|_{[t+h,T]},\qquad \forall u_2\in\mcU_{t+h}.
\end{align*}
We let $(\mcO_i)_{i\geq 1}\subset\mcB(\R^n)$ be a partition of $\R^n$ such that $(1+\sup_{x\in\mcO_i}|x|^\rho){\rm diam}(\mcO_i)\leq\epsilon$, then by Proposition~\ref{prop:W-k-cont} there is a $C> 0$ such that $|\Vlow(t+h,x)-\Vlow(t+h,x')|\leq C\epsilon$ for all $i\geq 1$ and $x,x'\in\mcO_i$. We pick $x_i\in\mcO_i$ and have by the same pasting property as above that there is for each $i\geq 1$ and $j\in \{0,\ldots,k\}$, a $u^\epsilon_{2,i,j}\in\mcU_{t+h}^j$ such that
\begin{align*}
  \Vlow^{j}(t+h,x_i) &\leq J(t+h,x_i,u^\epsilon_{2,i,j},\alpha^S_2(u^\epsilon_{2,i,j}))+\epsilon.
\end{align*}
Consequently,
\begin{align*}
&\Vlow^{k-N^\epsilon_1}(t+h,X^{t,x;u^\epsilon_1,\alpha_1^S(u^\epsilon_1)}_{t+h}) \leq \sum_{i\geq 1}\ett_{[X^{t,x;u^\epsilon_1,\alpha_1^S(u^\epsilon_1)}_{t+h}\in \mcO_i]}\Vlow^{k-N^\epsilon_1}(t+h,x_i)+C\epsilon
\\
&\leq\sum_{i\geq 1}\sum_{j=0}^k\ett_{[k-N^\epsilon_1=j]}\ett_{[X^{t,x;u^\epsilon_1,\alpha_1^S(u^\epsilon_1)}_{t+h}\in \mcO_i]}J(t+h,x_i,u^\epsilon_{2,i,j},\alpha^S_2(u^\epsilon_{2,i,j}))+C\epsilon
\\
&\leq\sum_{i\geq 1}\sum_{j=0}^k\ett_{[k-N^\epsilon_1=j]}\ett_{[X^{t,x;u^\epsilon_1,\alpha_1^S(u^\epsilon_1)}_{t+h}\in \mcO_i]}J(t+h,X^{t,x;u_1^\epsilon,\alpha_1^S(u_1^\epsilon)}_{t+h},u^\epsilon_{2,i,j},\alpha^S(u^\epsilon))+C\epsilon,
\end{align*}
with
\begin{align*}
  u^\epsilon:=u_1^\epsilon\oplus\sum_{i\geq 1}\sum_{j=0}^k\ett_{[k-N^\epsilon_1=j]}\ett_{X^{t,x;u_1^\epsilon,\alpha_1^S(u_1^\epsilon)}_{t+h}\in\mcO_i]}u_{2,i,j}.
\end{align*}
Using first comparison and then the stability property for BSDEs we find that
\begin{align*}
\Vlowh^{k}(t,x) &\leq G_{t,t+h}^{t,x;u^\epsilon_1,\alpha_1^S(u^\epsilon_1)}[\sum_{i\geq 1}\sum_{j=0}^k\ett_{[k-N^\epsilon_1=j]}\ett_{[X^{t,x;u^\epsilon_1,\alpha_1^S(u^\epsilon_1)}_{t+h}\in \mcO_i]}J(t+h,X^{t,x;u_1^\epsilon,\alpha_1^S(u_1^\epsilon)}_{t+h},u^\epsilon_{2,i,j},\alpha^S(u^\epsilon))+C\epsilon]
\\
&\leq J(t,x;u^\epsilon,\alpha^S(u^\epsilon))+C\epsilon
\\
&\leq \esssup_{u\in\mcU^k}J(t,x;u^\epsilon,\alpha^S(u))+C\epsilon,
\end{align*}
where $C>0$ only depends on the coefficients of the BSDE. Now, as this holds for all $\alpha^S\in \mcA^S_t$ we conclude that $\Vlowh^k(t,x)\leq \Vlow^k(t,x)+C\epsilon$, but $\epsilon>0$ was arbitrary and the result follows.\qed\\

The opposite inequality and its proof are classical (see \eg Proposition 1.10 in \cite{FlemSoug89} and Proposition 3.1 in \cite{TangHouSWgame}) and we give the proof only for the sake of completeness.

\begin{lem}\label{lem:W-lessthan-Wh}
$\Vlow^k(t,x)\leq \Vlowh^k(t,x)$.
\end{lem}

\noindent\emph{Proof.} We again fix an $\epsilon>0$ and let $(\mcO_i)_{i\geq 1}$ be defined as above. We pick an $x_i\in\mcO_i$ for each $i\geq 1$ and note that there is a $\alpha_{2,i,j}^S\in\mcA^S_{t+h,T}$ (see \cite{Buckdahn08} Lemma 4.5) such that
\begin{align*}
  \Vlow^j(t+h,x_i)\geq J(t+h,x_i;u_2,\alpha^S_{2,i,j}(u_2))-\epsilon,
\end{align*}
for all $u_2\in\mcU^j_{t+h}$. Moreover, there is an $\alpha^S_1\in\mcA_{t,t+h}^S$ such that
\begin{align*}
  \Vlowh^k(t,x)\geq G_{t,t+h}^{t,x;u_1,\alpha^S_1(u_1)}[\Vlow^{k-N_1}(t+h,X^{t,x;u_1,\alpha^S_1(u_1)}_{t+h})]-\epsilon,
\end{align*}
for all $u_1\in\mcU^k_{t,t+h}$, where $N_1$ is the number of interventions in $u_1$. Now, each $u=(\tau_i,\beta_i)_{1\leq i\leq N}\in\mcU^k_t$ can be uniquely decomposed as $u=u_1\oplus u_2$ with $u_1\in\mcU^k_{t,t+h}$ (with $N_1:=\max\{j\geq 0:\tau_j\leq t+h\}$ interventions) and $u_2\in\mcU^k_{t+h}$ (with first intervention at $\tau^2_1>t+h$). Then,
\begin{align*}
  \Vlowh^k(t,x)&\geq G_{t,t+h}^{t,x;u_1,\alpha^S_1(u_1)}[\Vlow^{k-N_1}(t+h,X^{t,x;u_1,\alpha^S_1(u_1)}_{t+h})]-\epsilon
  \\
  &=G_{t,t+h}^{t,x;u_1,\alpha^S_1(u_1)}[\sum_{j=0}^k\ett_{[k-N_1=j]} \Vlow^{j}(t+h,X^{t,x;u_1,\alpha^S_1(u_1)}_{t+h})]-\epsilon
  \\
  &\geq G_{t,t+h}^{t,x;u_1,\alpha^S_1(u_1)}[\sum_{j=0}^k\ett_{[k-N_1=j]}\sum_{i\geq 1}\ett_{[X^{t,x;u_1,\alpha^S_1(u_1)}_{t+h}\in\mcO_i]} \Vlow^{j}(t+h,x_i)]-C\epsilon
  \\
  &\geq G_{t,t+h}^{t,x;u_1,\alpha^S_1(u_1)}[\sum_{j=0}^k\ett_{[k-N_1=j]}\sum_{i\geq 1}\ett_{[X^{t,x;u_1,\alpha^S_1(u_1)}_{t+h}\in\mcO_i]} J(t+h,x_i;u_2,\alpha^S_{2,i,j}(u_2))]-C\epsilon
  \\
  &\geq G_{t,t+h}^{t,x;u_1,\alpha^S_1(u_1)}[\sum_{j=0}^k\ett_{[k-N_1=j]}\sum_{i\geq 1}\ett_{[X^{t,x;u_1,\alpha^S_1(u_1)}_{t+h}\in\mcO_i]} J(t+h,X^{t,x;u_1,\alpha^S(u_1)}_{t+h};u_2,\alpha^S_{2,i,j}(u))]-C\epsilon
  \\
  &=J(t,x;u,\alpha_1^S(u_1)\oplus_{t+h}\alpha_2^S(u))-C\epsilon,
\end{align*}
with
\begin{align*}
  \alpha_2^S(u):=\sum_{j=0}^k\ett_{[k-N_1=j]}\sum_{i\geq 1}\ett_{[X^{t,x;u_1,\alpha^S_1(u_1)}_{t+h}\in\mcO_i]}\alpha^S_{2,i,j}(u_2).
\end{align*}
Since $u\mapsto \alpha^S(u):=u\mapsto \alpha_1^S(u_1)\oplus_{t+h}\alpha_2^S(u)\in\mcA^S_{t}$, we conclude that $\Vlowh^k(t,x)\geq \Vlow(t,x)-C\epsilon$, where $C>0$ does not depend on $\epsilon>0$ which in turn was arbitrary and the result follows.\qed\\

Similarly, letting $\Vuph^k$ denote the right-hand-side of \eqref{ekv:dynp-U-trunk}, we find $\Vuph^k(t,x)=\Vup^k(t,x)$ for each $(t,x)\in[0,T]\times \R^n$ and the statement in Proposition~\ref{prop:dynp-trunk} follows.\qed\\

%Having shown continuity and that a dynamic programming principle holds for the truncated value functions we move on to consider the unlimited setting.

\subsection{A DPP for the general case}
We turn now to the general case where there is no restriction on the number of interventions in the impulse control. Before taking the limit as $k\to\infty$ in $\Vlow^k$ and $\Vup^k$, we need to delimit the set of impulse controls:
\begin{defn}\label{def:sens-controls}
For $(t,x)\in [0,T]\times \R^n$ and $\alpha^S\in\mcA_t^S$ we let $\bar\mcU_{t,x,\alpha^S}$ be the set of all $u\in\mcU_t$ such that $Y^{t,x;u,\alpha^S(u)}_s\geq Y^{t,x;u_{-s},\alpha^S(u_{-s})}_s$, $\Prob$-a.s., for all $s\in[t,T]$.

Moreover, we let $\bar\mcU^{S,k}_{t,x}$ be the subset of all $u^S\in\mcU^{S,k}_t$ such that for each $\alpha\in\mcA_t$ and $s\in [t,T]$,
\begin{align*}
  \essinf_{\tilde\alpha\in\mcA_s}Y^{t,x;u^S(\alpha\oplus_s\tilde\alpha),\alpha\oplus_s\tilde\alpha}_s\geq \essinf_{\tilde\alpha\in\mcA_s}Y^{t,x;(u^S(\alpha))_{s-},\alpha\oplus_s\tilde\alpha}_s
\end{align*}
$\Prob$-a.s.
\end{defn}
Given an $\alpha^S\in\mcA_t^S$ we note that the set $\bar\mcU_{t,x,\alpha^S}$ consists of all controls $u$ where it is never (on average) beneficial to abandon $u$ and stop intervening on the system for the remainder of the period. Similarly, $\bar\mcU^S_{t,x}$ is the set of strategies where, given that the opponent acts rationally, it will never be beneficial to abandon $u$ and stop intervening. The usefulness of the above definitions in our case lies in the fact that they allow us to bound the corresponding solution to \eqref{ekv:non-ref-bsde} from below by an expression that does not involve intervention costs. In particular, we have whenever $\alpha^S\in\mcA_t$ and $u\in \bar\mcU_{t,x,\alpha^S}$, that
\begin{align}\label{ekv:sens-useful-C}
Y^{t,x;u,\alpha^S(u)}_s\geq Y^{t,x;u_{-s},\alpha^S(u_{-s})}_s\geq\essinf_{\tilde\alpha\in\mcA_s}Y^{t,x;u_{-s},\alpha^S(u_{-s})\oplus_s\tilde\alpha}_s
\end{align}
for all $s\in [t,T]$, and similarly when $u^S\in\bar\mcU^{S,k}_{t,x}$ we have
\begin{align}\label{ekv:sens-useful-S}
Y^{t,x;u^S(\alpha),\alpha}_s\geq \essinf_{\tilde\alpha\in\mcA_s}Y^{t,x;u^S(\alpha\oplus_s\tilde\alpha),\alpha\oplus_s\tilde\alpha}_s\geq \essinf_{\tilde\alpha\in\mcA_s}Y^{t,x;(u^S(\alpha))_{s-},\alpha\oplus_s\tilde\alpha}_s
\end{align}
for all $\alpha\in\mcA_t$ and $s\in [t,T]$.

The following lemma shows that these sets contain all relevant impulse controls and strategies, respectively.
\begin{lem}\label{lem:sens-is-opt}
We have
\begin{align*}
  \Vlow(t,x)=\essinf_{\alpha^S\in\mcA^S_t}\esssup_{u\in\bar\mcU_{t,x,\alpha^S}}J(t,x;u,\alpha^S(u))
\end{align*}
and
\begin{align*}
  \Vup(t,x)=\esssup_{u^S\in\bar\mcU^{S}_{t,x}}\essinf_{\alpha\in\mcA_t}J(t,x;u^S(\alpha),\alpha).
\end{align*}
\end{lem}

\noindent\emph{Proof.} For any $\alpha^S\in\mcA^S_t$ and arbitrary $u\in\mcU_t\setminus \bar\mcU_{t,x,\alpha^S}$ we let
\begin{align*}
  \chi:=\inf\big\{s\geq t: Y^{t,x;u,\alpha^S(u)}_s\leq Y^{t,x;u_{s-},\alpha^S(u_{s-})}_s\big\}\wedge T.
\end{align*}
Assumption~\ref{ass:on-coeff}.\ref{ass:on-coeff-@end} implies that $Y^{t,x;u,\alpha^S(u)}_T\leq Y^{t,x;u_{T-},\alpha^S(u_{T-})}_T$ and we get that with
\begin{align*}
B_1:=\{\omega:Y^{t,x;u,\alpha^S(u)}_\chi\leq Y^{t,x;u_{\chi-},\alpha^S(u_{\chi-})}_\chi\}\in\mcF_\chi
\end{align*}
and
\begin{align*}
B_2:=\{\omega:Y^{t,x;u,\alpha^S(u)}_{\chi+}\leq Y^{t,x;u_{\chi},\alpha^S(u_{\chi})}_{\chi+}\}\cap B_1^c\in\mcF_\chi,
\end{align*}
the set $(B_1\cup B_2)^c$ is $\Prob$-negligible.

Moreover, since
\begin{align*}
Y^{t,x;u,\alpha^S(u)}_{\chi+}-Y^{t,x;u,\alpha^S(u)}_{\chi}=Y^{t,x;u_{\chi},\alpha^S(u_{\chi})}_{\chi+}-Y^{t,x;u_{\chi},\alpha^S(u_{\chi})}_{\chi},
\end{align*}
it follows that on $B_2$ we have $Y^{t,x;u,\alpha^S(u)}_{\chi}\leq Y^{t,x;u_{\chi},\alpha^S(u_{\chi})}_{\chi}$ and we conclude that letting $\tilde u:=\ett_{B_1}u_{\chi-}+\ett_{B_2}u_{\chi}$ we have $Y^{t,x;u,\alpha^S(u)}_\chi\leq Y^{t,x;\tilde u,\alpha^S(\tilde u)}_\chi$ $\Prob$-a.s. By comparison we thus find that $Y^{t,x;u,\alpha^S(u)}_s \leq Y^{t,x;\tilde u,\alpha^S(\tilde u)}_s$, $\Prob$-a.s.~for all $s\in [t,\chi]$. In particular, this gives that $\tilde u\in \bar\mcU_{t,x,\alpha^S}$ and $Y^{t,x;u,\alpha^S(u)}_t \leq Y^{t,x;\tilde u,\alpha^S(\tilde u)}_t$ from which we conclude that any $u\in\mcU_t\setminus \bar\mcU_{t,x,\alpha^S}$ is dominated by an element of $\bar\mcU_{t,x,\alpha^S}$. Since this holds for any $\alpha^S\in\mcA^S_t$, we have that
\begin{align*}
  \essinf_{\alpha^S\in\mcA^S_t}\esssup_{u\in\mcU_{t}}J(t,x;u,\alpha^S(u))= \essinf_{\alpha^S\in\mcA^S_t}\esssup_{u\in\bar\mcU_{t,x,\alpha^S}}J(t,x;u,\alpha^S(u)),
\end{align*}
proving the first statement.

For the second statement we fix $u^S\in\mcU_t^{S}\setminus \bar\mcU_{t,x}^{S}$ and $\alpha\in\mcA_t$. We then set $u=(\tau_i,\beta_i)_{1\leq i\leq N}:=u^S(\alpha)$ and let
\begin{align*}
  N(\alpha):=\min\Big\{j\geq 0: \essinf_{\tilde\alpha\in\mcA_{\tau_j}}Y^{t,x;u^S(\alpha\oplus_{\tau_j}\tilde\alpha),\alpha\oplus_{\tau_j}\tilde\alpha}_{\tau_j}\leq \essinf_{\tilde\alpha\in\mcA_{\tau_j}}Y^{t,x;[u^S(\alpha)]_{j-1},\alpha\oplus_{\tau_j}\tilde\alpha}_{\tau_j}\Big\}.
\end{align*}
Furthermore, we define $\tilde u^S\in\mcU^S_t$ as $\tilde u^S(\alpha):=[u^S(\alpha)]_{N(\alpha)-1}$ and let $\chi(\alpha):=\tau_{N(\alpha)}\wedge T$. By definition we have
\begin{align}\label{ekv:end-rel}
  \essinf_{\tilde\alpha\in\mcA_{\chi(\alpha)}} Y^{t,x;u^S(\alpha\oplus_{\chi(\alpha)}\tilde\alpha),\alpha\oplus_{\chi(\alpha)}\tilde\alpha}_{\chi(\alpha)}\leq \essinf_{\tilde\alpha\in\mcA_{\chi(\alpha)}} Y^{t,x;\tilde u^S(\alpha),\alpha\oplus_{\chi(\alpha)}\tilde\alpha}_{\chi(\alpha)}.
\end{align}
For $\epsilon>0$ and $s\in[t,T]$ we let $\mcA^{\epsilon}_{s}$ be the subset of all $\hat\alpha\in\mcA_t$ with $\hat\alpha=\alpha$ on $[t,s)$ such that
\begin{align*}
  \essinf_{\tilde\alpha\in\mcA_{\chi(\hat\alpha)\vee s}}Y^{t,x;u^S(\hat\alpha\oplus_{\chi(\hat\alpha)\vee s}\tilde\alpha),\hat\alpha\oplus_{\chi(\hat\alpha)\vee s}\tilde\alpha}_{\chi(\hat\alpha)\vee s}\geq Y^{t,x;u^S(\hat\alpha),\hat\alpha}_{\chi(\hat\alpha)\vee s}-\epsilon
\end{align*}
and similarly let $\tilde\mcA^{\epsilon}_{s}$ be the subset of all $\hat\alpha\in\mcA_t$ with $\hat\alpha=\alpha$ on $[t,s)$ such that
\begin{align*}
  \essinf_{\tilde\alpha\in\mcA_{\chi(\hat\alpha)\vee s}}Y^{t,x;\tilde u^S(\hat\alpha\oplus_{\chi(\hat\alpha)\vee s}\tilde\alpha),\hat\alpha\oplus_{\chi(\hat\alpha)\vee s}\tilde\alpha}_{\chi(\hat\alpha)\vee s}\geq Y^{t,x;\tilde u^S(\hat\alpha),\hat\alpha}_{\chi(\hat\alpha)\vee s}-\epsilon.
\end{align*}
Then, we can repeat the arguments in Lemma~\ref{lem:Wh-lessthan-W} to conclude that for all $s\in [t,T]$, the sets $\mcA^{\epsilon}_s$ and $\tilde\mcA^{\epsilon}_s$ are non-empty and comparison implies that
\begin{align}\label{ekv:A-eps-opt}
 \essinf_{\tilde\alpha\in\mcA_s}Y^{t,x; u^S(\alpha\oplus_s\tilde\alpha),\alpha\oplus_s\tilde\alpha}_s = \essinf_{\hat\alpha\in\mcA^{\epsilon}_s}Y^{t,x; u^S(\hat\alpha),\hat\alpha}_s
\end{align}
and
\begin{align}\label{ekv:tilde-A-eps-opt}
  \essinf_{\hat\alpha\in\mcA_s}Y^{t,x;\tilde u^S(\alpha\oplus_s\hat\alpha),\alpha\oplus_s\hat\alpha}_s=\essinf_{\hat\alpha\in\tilde\mcA^{\epsilon}_s}Y^{t,x; \tilde u^S(\hat\alpha),\hat\alpha}_s.
\end{align}
Moreover, for $\hat\alpha\in\mcA^{\epsilon}_s$ and $\tilde\alpha\in\tilde\mcA^{\epsilon}_s$ with $\hat\alpha=\tilde \alpha$ on $[t,\chi(\hat\alpha))$ we have by \eqref{ekv:end-rel}, that
\begin{align*}
 Y^{t,x;u^S(\hat\alpha),\hat\alpha}_{\chi(\hat\alpha)} \leq Y^{t,x;\tilde u^S(\tilde\alpha),\tilde\alpha}_{\chi(\tilde\alpha)}+\epsilon.
\end{align*}
and using comparison together with stability implies that
\begin{align*}
  \essinf_{\hat\alpha\in\mcA^\epsilon_s}Y^{t,x;u^S(\hat\alpha),\hat\alpha}_s\leq \essinf_{\hat\alpha\in\tilde\mcA^\epsilon_s}Y^{t,x;\tilde u^S(\hat\alpha),\hat\alpha}_s+C\epsilon
\end{align*}
for all $s\in [t,\chi(\alpha)]$. In particular, since $\epsilon>0$ was arbitrary, letting $s=t$ and using \eqref{ekv:A-eps-opt} and \eqref{ekv:tilde-A-eps-opt} gives that
\begin{align*}
   \essinf_{\alpha\in\mcA_t}J(t,x;u^S(\alpha),\alpha)\leq \essinf_{\alpha\in\mcA_t} J(t,x;\tilde u^S(\alpha),\alpha)
\end{align*}
and we conclude that $\tilde u^S$ dominates $u^S$. On the other hand, by a similar argument we find that
\begin{align*}
  \essinf_{\tilde\alpha\in\mcA_s}Y^{t,x;\tilde u^S(\alpha\oplus_s\tilde\alpha),\alpha\oplus_s\tilde\alpha}_s\geq \essinf_{\tilde\alpha\in\mcA_s}Y^{t,x;(u^S(\alpha))_{s-},\alpha\oplus_s\tilde\alpha}_s =\essinf_{\tilde\alpha\in\mcA_s}Y^{t,x;(\tilde u^S(\alpha))_{s-},\alpha\oplus_s\tilde\alpha}_s
\end{align*}
for all $s\in[t,\chi(\alpha)]$ and since $(\tilde u^S(\alpha))_{s-}=\tilde u^S(\alpha)$ on $(\chi(\alpha),T]$ we conclude that
$\tilde u^S\in\bar\mcU^S_{t,x}$ and the assertion follows.\qed\\

In particular, we may w.l.o.g.~restrict our attention to impulse controls (resp. strategies) in Definition~\ref{def:sens-controls}. The following result relates the number of interventions in these impulse controls and strategies to the magnitude of the initial value and is central in deriving continuity of $\Vlow$ and $\Vup$.

\begin{lem}\label{lem:EN-bound}
There is a constant $C>0$ such that
\begin{align}\label{ekv:N-bound}
  \E[N]\leq C(1+|x|^\rho)
\end{align}
for all $\alpha^S\in \mcA^S_t$ and $u\in\bar\mcU_{t,x,\alpha^S}$. Moreover, \eqref{ekv:N-bound} also holds for $u=u^S(\alpha)$ whenever $u^S\in\bar\mcU^{S}_{t,x}$ and $\alpha\in\mcA_t$.
\end{lem}

\noindent \emph{Proof.} Both statements will follow by a similar argument and we set $\alpha:=\alpha^S(u)$ (resp. $u:=u^S(\alpha)$). To simplify notation we let $(X,Y,Z):=(X^{t,x;u,\alpha},Y^{t,x;u,\alpha},Z^{t,x;u,\alpha})$ and $X^j=X^{t,x;[u]_j,\alpha}$ and get that
\begin{align*}
Y_s&=\psi(X_T)+\int_s^{T}f(r,X_r,Y_r,Z_r,\alpha_r)dr-\int_s^{T} Z_rdW_r-\sum_{\tau_j\geq s}\ell(\tau_j,X^{j-1}_{\tau_j},\beta_j).
\end{align*}
Letting
\begin{align*}
 \zeta^{u,\alpha}_1(s):=\frac{f(s,X^{t,x;u,\alpha}_s,Y^{t,x;u,\alpha}_s,Z^{t,x;u,\alpha}_s,\alpha_s) -f(s,X^{t,x;u,\alpha}_s,0,Z^{t,x;u,\alpha}_s,\alpha_s)}{Y^{t,x;u,\alpha}_s}\ett_{[Y_s\neq 0]}
\end{align*}
and
\begin{align*}
 \zeta^{u,\alpha}_2(s):=\frac{f(s,X^{t,x;u,\alpha}_s,0,Z^{t,x;u,\alpha}_s,\alpha_s) - f(s,X^{t,x;u,\alpha}_s,0,0,\alpha_s)}{|Z^{t,x;u,\alpha}_s|^2}(Z^{t,x;u,\alpha}_s)^\top
\end{align*}
we have by the Lipschitz continuity of $f$ that $|\zeta^{u,\alpha}_1(s)|\vee|\zeta^{u,\alpha}_2(s)|\leq k_f$. Using Ito's formula we find that
\begin{align*}
Y_s&=R^{u,\alpha}_{s,T}\psi(X_T)+\int_s^{T}R^{u,\alpha}_{s,r}f(r,X_r,0,0,\alpha_r)dr -\int_s^TR^{u,\alpha}_{s,r}Z_rdW_r-\sum_{j=1}^N R^{u,\alpha}_{s,\tau_{j}}\ett_{[\tau_j\geq s]}\ell(\tau_j,X^{j-1}_{\tau_{j}},\beta_j)
\end{align*}
with $R^{u,\alpha}_{s,r}:=e^{\int_s^{r}(\zeta^{u,\alpha}_1(v)-\frac{1}{2}|\zeta^{u,\alpha}_2(v)|^2)dv+\frac{1}{2}\int_s^{r}\zeta^{u,\alpha}_2(v)dW_v}$. Since the intervention costs are positive, taking the conditional expectation on both sides gives
\begin{align*}
Y_s&\leq\E\Big[R^{u,\alpha}_{s,T}\psi(X_T)+\int_s^{T}R^{u,\alpha}_{s,r}f(r,X_r,0,0,\alpha_r)dr \Big|\mcF_s\Big]
\\
&\leq C(1+|X_s|^\rho)
\end{align*}
On the other hand, by \eqref{ekv:sens-useful-C} (resp.~\eqref{ekv:sens-useful-S}) we have
\begin{align*}
Y_s&\geq \essinf_{\tilde\alpha\in\mcA_s}\E\Big[R^{u_{s-},\alpha\oplus_s\tilde\alpha}_{s,T}\psi(X^{t,x;u_{s-},\alpha\oplus_s\tilde\alpha}_T) +\int_s^{T}R^{u_{s-},\alpha\oplus_s\tilde\alpha}_{s,r}f(r,X^{t,x;u_{s-},\alpha\oplus_s\tilde\alpha}_r,0,0,\alpha_r)dr \Big|\mcF_s\Big]
\\
&\geq -C(1+|X^{t,x;u_{s-},\alpha}_s|^\rho).
\end{align*}
Proposition~\ref{prop:SDEmoment} then gives
\begin{align*}
\E\Big[\sup_{s\in[t,T]}|Y_s|^2\Big]&\leq C(1+|x|^{2\rho}).
\end{align*}
Next, we derive a bound on the $\mcH^2$-norm of $Z$. Applying Ito's formula to $|Y_s|^{2}$ we get
\begin{align}\nonumber
|Y_t|^{2}+\int_t^T| Z_s|^2ds&=\psi^2(X_T)+\int_t^TY_sf(s,X_s,Y_s,Z_s,\alpha_s)ds-2\int_t^T Y_sZ_sdW_s
\\
&\quad-\sum_{j=1}^N(2Y^{j-1}_{\tau_{j}}\ell(\tau_j,X^{j-1}_{\tau_{j}},\beta_j)+\ell^2(\tau_j,X^{j-1}_{\tau_{j}},\beta_j)),\label{ekv:from-ito}
\end{align}
where $Y^{j-1}$ is $Y$ without the $j-1$ first intervention costs. Since the intervention costs are nonnegative, we have
\begin{align*}
-\sum_{j=1}^N(2Y^{j-1}_{\tau_{j}}\ell(\tau_j,X^{j-1}_{\tau_{j}},\beta_j)+\ell^2(\tau_j,X^{j-1}_{\tau_{j}},\beta_j)) &\leq 2\sup_{s\in [t,T]}|Y_{s}|\sum_{j=1}^N\ell(\tau_j,X^{j-1}_{\tau_{j}},\beta_j)
\\
&\leq \kappa \sup_{s\in [t,T]}|Y_{s}|^2+\frac{1}{\kappa}\Big(\sum_{j=1}^N\ell(\tau_j,X^{j-1}_{\tau_{j}},\beta_j)\Big)^2
\end{align*}
for any $\kappa>0$. Inserted in \eqref{ekv:from-ito} and using the Lipschitz property of $f$ this gives
\begin{align*}
|Y_t|^{2}+\int_t^T| Z_s|^2ds&\leq \psi^2(X_T)+(C+\kappa)\sup_{s\in[t,T]}|Y_s|^2+\int_t^T(|f(s,X_s,0,0,\alpha_s)|^2+\frac{1}{2}|Z_s|^2)ds
\\
&\quad -2\int_t^T Y_sZ_sdW_s+\frac{1}{\kappa}\Big(\sum_{j=1}^N\ell(\tau_j,X^{j-1}_{\tau_{j}},\beta_j)\Big)^2.
\end{align*}
Now, as $u\in\mcU$, it follows that the stochastic integral is uniformly integrable and thus a martingale. To see this, note that the Burkholder-Davis-Gundy inequality gives
\begin{align*}
\E\Big[\sup_{s\in[t,T]}\Big|\int_t^s Y_rZ_rdW_r\Big|\Big]\leq C\E\Big[\Big(\int_t^T |Y_sZ_s|^2ds\Big)^{1/2}\Big]\leq C\E\Big[\sup_{s\in [t,T]}|Y_s|^2+\int_t^T |Z_s|^2ds\Big].
\end{align*}
Taking expectations on both sides thus gives
\begin{align*}
\E\Big[\int_t^T| Z_s|^2ds\Big]&\leq C(1+\kappa)(1+|x|^{2\rho})+\frac{2}{\kappa}\E\Big[\Big(\sum_{j=1}^N\ell(\tau_j,X^{j-1}_{\tau_{j}},\beta_j)\Big)^2\Big].
\end{align*}
Finally,
\begin{align*}
\E[N]&\leq \frac{1}{\delta}\E\Big[\Big(\sum_{j=1}^N\ell(\tau_j,X^{j-1}_{\tau_{j}},\beta_j)\Big)^2\Big]^{1/2}
\end{align*}
and
\begin{align*}
\E\Big[\Big(\sum_{j=1}^N\ell(\tau_j,X^{j-1}_{\tau_{j}},\beta_j)\Big)^2\Big]&\leq C\E\Big[|Y_t|^2+|\psi(X_T)|^2+\int_t^{T}|f(r,X_r,Y_r,Z_r,\alpha_r)|^2dr+\int_t^{T}|Z_r|^2dr\Big]
\\
&\leq C\E\Big[|\psi(X_T)|^2+\sup_{s\in[t,T]}|Y_s|^2+\int_t^{T}(|f(s,X_s,0,0,\alpha_s)|^2+|Z_s|^2)ds\Big]
\\
&\leq C(1+\kappa)(1+|x|^{2\rho})+\frac{C}{\kappa}\E\Big[\Big(\sum_{j=1}^N\ell(\tau_j,X^{j-1}_{\tau_{j}},\beta_j)\Big)^2\Big]
\end{align*}
from which \eqref{ekv:N-bound} follows by choosing $\kappa$ sufficiently large.\qed\\

\begin{lem}\label{lem:unif-conv}
There is a $C>0$ such that for all $k\geq 1$ we have
\begin{align*}
  \Vlow(t,x)-\Vlow^{k}(t,x)\leq \frac{C(1+|x|^{2\rho})}{\sqrt{k}}
\end{align*}
for all $(t,x)\in[0,T]\times \R^n$. In particular, the sequence $\{\Vlow^k\}_{k\geq 0}$ converges to $\Vlow$, uniformly on compact subsets of $[0,T]\times\R^n$.
\end{lem}

\noindent\emph{Proof.} For each $\alpha^S\in\mcA^S$ and $\epsilon>0$ there is by Lemma~\ref{lem:sens-is-opt} a $u^{\epsilon}=(\tau_i^{\epsilon}, \beta_i^{\epsilon})_{1\leq i\leq N^{\epsilon}}\in\bar\mcU_{t,x,\alpha^S}$ such that
\begin{align}\label{ekv:u-eps-optimal}
  \esssup_{u\in\mcU_t}J(t,x;u,\alpha^S(u))\leq J(t,x;u^{\epsilon},\alpha^S(u^{\epsilon}))+\epsilon/2,
\end{align}
$\Prob$-a.s.
%and let $(\mcY,\mcZ)$ solve the BSDE
%\begin{align*}
%\mcY_s&=\psi(X^{t,x;u^{\epsilon},\alpha}_T)+\int_s^{T}f(r,X^{t,x;u^{\epsilon},\alpha}_r,\mcY_r,\mcZ_r)dr-\int_s^{T} \mcZ_rdW_r -\Xi^{t,x;u^{\epsilon},\alpha}_{T+} +\Xi^{t,x;u^{\epsilon},\alpha}_s.
%\end{align*}
Now, let $(Y,Z)=(Y^{{t,x;u^{\epsilon},\alpha^S(u^{\epsilon})}},Z^{{t,x;u^{\epsilon},\alpha^S(u^{\epsilon})}})$ and for $k\geq 0$, set\\ $(\hat Y,\hat Z)=(Y^{{t,x;[u^{\epsilon}]_{k},\alpha^S([u^{\epsilon}]_{k})}},Y^{{t,x;[u^{\epsilon}]_{k},\alpha^S([u^{\epsilon}]_{k})}})$, %solve
%\begin{align*}
%\hat\mcY_s&=\psi(X^{t,x;[u^{\epsilon}]_{k},\alpha}_T)+\int_s^{T}f(r,X^{t,x;[u^{\epsilon}]_{k},\alpha}_r,\hat\mcY_r,\hat\mcZ_r)d r-\int_s^{T} \hat\mcZ_rdW_r-\Xi^{t,x;[u^{\epsilon}]_{k},\alpha}_{T+} +\Xi^{t,x;[u^{\epsilon}]_{k},\alpha}_s,
%\end{align*}
where we recall that $[u]_l$ is the truncation of $u$ to the first $l$ interventions. %Then
%\begin{align*}
%  0 \leq \Vlow^{k}(t,x)-\Vlow^{k}(t,x)\leq \E[Y_t-\hat Y_t]+\epsilon.
%\end{align*}
As $\alpha^S([u^{\epsilon}]_{k})=\alpha^S(u^{\epsilon})$ on $[0,\tau^{\epsilon}_{k+1})\cap[0,T]$ we have, with $X:=X^{t,x;u^{\epsilon},\alpha^S(u^{\epsilon})}$ and $\hat X:=X^{t,x;[u^{\epsilon}]_{k},\alpha^S([u^{\epsilon}]_{k})}$, that $\hat X_s=X_s$ for all $s\in[0,\tau^{\epsilon}_{k+1})\cap[0,T]$. Letting $\alpha:=\alpha^S(u^{\epsilon})$ and $\hat\alpha:=\alpha^S([u^{\epsilon}]_{k})$, this gives
\begin{align*}
Y_t-\hat Y_t&=\psi(X_T)-\psi(\hat X_T) +\int_t^{T}(f(s,X_s,Y_s,Z_s,\alpha_s)- f(s,\hat X_s,\hat Y_s,\hat Z_s,\hat\alpha_s))ds
\\
&\quad -\int_t^{T}(Z_s- \hat Z_s)dW_s+\Xi^{t,x;[u^{\epsilon}]_{k},\alpha}_{T+}-\Xi^{t,x;u^{\epsilon},\alpha^S(u^{\epsilon})}_{T+}
\\
&\leq \ett_{[N^{\epsilon}>k]}\Big(R_T(\psi(X_T)-\psi(\hat X_T)) +\int_t^{T}R_s(f(s,X_s,Y_s,Z_s,\alpha_s)- f(s,\hat X_s,Y_s,Z_s,\hat\alpha_s))ds\Big)
\\
&\quad-\int_t^{T}R_s(Z_s- \hat Z_s)dW_s
\end{align*}
for some $R_{s}:=e^{\int_t^{s}(\zeta_1(r)-\frac{1}{2}|\zeta_2(r)|^2)dr+\frac{1}{2}\int_t^{s}\zeta_2(r)dW_r}$, with $|\zeta_1(r)|\vee|\zeta_2(r)|\leq k_f$. Taking expectation on both sides and using the Cauchy-Schwartz inequality gives
\begin{align*}
\E[Y_t-\hat Y_t]&\leq\E\Big[\ett_{[N^{\epsilon}>k]}\Big(R_T(\psi(X_T) - \psi(\hat X_T))
+\int_t^{T}R_s(f(s,X_s,Y_s,Z_s,\alpha_s)- f(s,\hat X_s,Y_s,Z_s,\hat\alpha_s))ds\Big) \Big]
\\
&\leq C(1+|x|^{\rho})\E\big[\ett_{[N^{\epsilon}>k]}\big]^{1/2}.
\end{align*}
Now, as $u^{\epsilon}\in\bar\mcU_{t,x,\alpha^S}$, Lemma~\ref{lem:EN-bound} implies that
\begin{align*}
\E\big[\ett_{[N^{\epsilon}>k]}\big]\leq  \frac{C(1+|x|^\rho)}{k}.
\end{align*}
Since $\alpha^S$ was arbitrary we can pick $\alpha^S\in\mcA_t$ such that
\begin{align*}
  \Vlow^{k}(t,x)\geq \esssup_{u\in\mcU^{k}_t}J(t,x;u,\alpha^S(u))-\epsilon/2
\end{align*}
and we find that
\begin{align*}
  \Vlow(t,x)-\Vlow^{k}(t,x)&\leq \esssup_{u\in\mcU_t} J(t,x;u,\alpha^S(u))- \esssup_{u\in\mcU^{k}_t} J(t,x;u,\alpha^S(u))+\epsilon/2
  \\
  &\leq J(t,x;u^{\epsilon},\alpha^S(u^{\epsilon}))-J(t,x;[u^{\epsilon}]_{k},\alpha^S([u^{\epsilon}]_{k}))+\epsilon
  \\
  &\leq \frac{C(1+|x|^{2\rho})}{\sqrt{k}}+\epsilon
\end{align*}
from which the desired inequality follows since $\epsilon>0$ was arbitrary. In particular, we find that $\Vlow^{k}$ converges uniformly on sets where $x$ is bounded.\qed\\

\begin{thm}\label{thm:W-is-cont}
$\Vlow$ is continuous and satisfies \eqref{ekv:dynp-W}
\end{thm}

\noindent\emph{Proof.} Since the sequence $\{\Vlow^k\}_{k\geq 0}$ is non-decreasing, Lemma~\ref{lem:unif-conv} implies that $\Vlow^k\nearrow\Vlow$ uniformly on compacts as $k\to\infty$. Hence, $\Vlow$ is continuous.

It remains to show that $\Vlow$ satisfies \eqref{ekv:dynp-W}. We have by \eqref{ekv:dynp-W-trunk} and comparison that
\begin{align*}
  \Vlow^k(t,x)&=\essinf_{\alpha^S\in\mcA^S_{t,t+h}}\esssup_{u\in\mcU^{k}_{t,t+h}}G_{t,t+h}^{t,x;u,\alpha^S(u)} [\Vlow^{k-N}(t+h,X^{t,x;u,\alpha^S(u)}_{t+h})]
  \\
  &\leq \essinf_{\alpha^S\in\mcA^S_{t,t+h}}\esssup_{u\in\mcU_{t,t+h}}G_{t,t+h}^{t,x;u,\alpha^S(u)}[\Vlow(t+h,X^{t,x;u,\alpha^S(u)}_{t+h})] =:\Vlowh(t,x)
\end{align*}
and it follows that $\Vlow\leq \Vlowh$. On the other hand, for each $\epsilon>0$ and any $\alpha^S\in\mcA^S$ we can repeat the argument in Lemma~\ref{lem:unif-conv} to find that there is a $k\geq 0$ such that% and a $u^\epsilon\in \mcU^{k}$ such that
\begin{align*}
\esssup_{u\in\mcU_{t,t+h}}G_{t,t+h}^{t,x;u,\alpha^S(u)}[\Vlow(t+h,X^{t,x;u,\alpha^S(u)}_{t+h})]%&\leq G_{t,t+h}^{t,x;u^\epsilon,\alpha^S(u^\epsilon)}[\barVlow(t+h,X^{t,x;u^\epsilon,\alpha^S(u^\epsilon)}_{t+h})]+\epsilon/2
%\\
&\leq \esssup_{u\in\mcU^{k}_{t,t+h}}G_{t,t+h}^{t,x;u,\alpha^S(u)}[\Vlow(t+h,X^{t,x;u,\alpha^S(u)}_{t+h})]+\epsilon/2.
\end{align*}
Moreover, for each $(\alpha,u)\in\mcA_t\times\mcU_t$, let $(\mcY^k,\mcZ^k)$ be the unique solution to
\begin{align*}
\mcY^k_s&=\Vlow^k(t+h,X^{t,x;u,\alpha}_{t+h})+\int_s^{t+h}f(r,X^{t,x;u,\alpha}_r,\mcY^k_r,\mcZ^k_r,\alpha_r)dr
\\
&\quad-\int_s^{t+h} \mcZ^k_rdW_r -\Xi^{t,x;u,\alpha}_{T+}+\Xi^{t,x;u,\alpha}_s
\end{align*}
while we assume that $(\mcY,\mcZ)$ satisfies
\begin{align*}
\mcY_s&=\Vlow(t+h,X^{t,x;u,\alpha}_{t+h})+\int_s^{t+h}f(r,X^{t,x;u,\alpha}_r,\mcY_r,\mcZ_r,\alpha_r)dr
\\
&\quad-\int_s^{t+h} \mcZ_rdW_r -\Xi^{t,x;u,\alpha}_{T+}+\Xi^{t,x;u,\alpha}_s.
\end{align*}
Then $\mcY_t\geq \mcY^k_t$ by comparison and
\begin{align}\nonumber
\mcY_t-\mcY^k_t&=\E\big[R_T(\Vlow(t+h,X^{t,x;u,\alpha}_{t+h}) - \Vlow^k(t+h,X^{t,x;u,\alpha}_{t+h}))\big]
\\
&\leq \frac{C}{\sqrt k}\E\big[R_T(1+|X^{t,x;u,\alpha}_{t+h}|^{2\rho})\big|\mcF_t\big],\label{ekv:G-stability}
\end{align}
with $R_{s}:=e^{\int_t^{s}(\zeta_1(r)-\frac{1}{2}|\zeta_2(r)|^2)dr+\frac{1}{2}\int_t^{s}\zeta_2(r)dW_r}$, where $|\zeta_1(r)|\vee|\zeta_2(r)|\leq k_f$. Since the right-hand side of the above inequality tends to 0 as $k\to\infty$ we conclude by taking the essential supremum over all $(\alpha^S,u)\in\mcA_t\times\mcU_t$ that there is a $k'\geq k$ such that
\begin{align*}
\esssup_{u\in\mcU^{k'}_{t,t+h}}G_{t,t+h}^{t,x;u,\alpha^S(u)}[\Vlow(t+h,X^{t,x;u,\alpha^S(u)}_{t+h})]\leq \esssup_{u\in\mcU^{k'}_{t,t+h}}G_{t,t+h}^{t,x;u,\alpha^S(u^\epsilon)}[\Vlow^{k'}(t+h,X^{t,x;u^\epsilon,\alpha^S(u^\epsilon)}_{t+h})]+\epsilon/2.
\end{align*}
for each $\alpha^S\in\mcA^S_t$. We conclude that $\Vlowh(t,x)\leq \Vlow^{k'}(t,x)+\epsilon\leq \Vlow(t,x)+\epsilon$ and since $\epsilon>0$ was arbitrary it follows that $\Vlow$ satisfies \eqref{ekv:dynp-W}.\qed\\

\begin{lem}\label{lem:unif-conv-U}
There is a constant $C>0$ such that for all $k\geq 1$ we have
\begin{align*}
  \Vup(t,x)-\Vup^{k'}(t,x)\leq \frac{C(1+|x|^{2\rho})}{\sqrt{k}}.
\end{align*}
for all $(t,x)\in[0,T]\times \R^n$. In particular, the sequence $\{\Vup^k\}_{k\geq 0}$ converges uniformly on compact subsets of $[0,T]\times\R^n$.
\end{lem}

\noindent\emph{Proof.} For each $\epsilon>0$ there is a $u^{\epsilon}=(\tau_i^{\epsilon},\beta_i^{\epsilon})_{1\leq i\leq N^{\epsilon}}\in\bar\mcU^{S}_{t,x}$ such that
\begin{align*}
  \esssup_{u^S\in\mcU^{S}_t}\essinf_{\alpha\in\mcA_t} J(t,x;u^S(\alpha),\alpha)\leq \essinf_{\alpha\in\mcA_t}J(t,x;u^{\epsilon}(\alpha),\alpha)+\epsilon,
\end{align*}
$\Prob$-a.s. Then, for $k\geq 0$,
\begin{align*}
  \Vup(t,x)-\Vup^{k}(t,x)&=\esssup_{u\in\mcU^{S}_{t}}\essinf_{\alpha\in\mcA_t} J(t,x;u(\alpha),\alpha)-\esssup_{u\in\mcU^{S,k}_t}\essinf_{\alpha\in\mcA_t} J(t,x;u(\alpha),\alpha)
  \\
  &\leq \essinf_{\alpha\in\mcA_t}J(t,x;u^{\epsilon}(\alpha),\alpha)-\essinf_{\alpha\in\mcA_t}J(t,x;[u^{\epsilon}(\alpha)]_{k},\alpha)+\epsilon
  \\
  & \leq \esssup_{\alpha\in\mcA_t}\{J(t,x;u^{\epsilon}(\alpha),\alpha)-J(t,x;[u^{\epsilon}(\alpha)]_{k'},\alpha)\}+\epsilon.
\end{align*}
By arguing as in the proof of Lemma~\ref{lem:unif-conv} the result now follows.\qed\\

\begin{thm}
$\Vup$ is continuous and satisfies \eqref{ekv:dynp-U}.
\end{thm}

%\noindent\emph{Proof.} The result follows by a similar argument to the one used in the proof of Theorem~\ref{thm:W-is-cont}.\qed\\

\noindent\emph{Proof.} As above we find that $\Vup^k\nearrow \Vup$ uniformly on compacts and conclude that $\Vup$ is continuous.

We have again by comparison that
\begin{align*}
  \Vup(t,x)\leq \esssup_{u^S\in\mcU^S_{t,t+h}}\essinf_{\alpha\in\mcA_{t,t+h}}G_{t,t+h}^{t,x;u^S(\alpha),\alpha}[\Vup(t+h,X^{t,x;u^S(\alpha),\alpha}_{t+h})] =:\Vuph(t,x).
\end{align*}
Moreover, for each $\epsilon>0$ there is a $k\geq 0$ such that
\begin{align*}
  \Vuph(t,x)\leq \esssup_{u^S\in\mcU^{S,k}_{t,t+h}}\essinf_{\alpha\in\mcA_{t,t+h}}G_{t,t+h}^{t,x;u^S(\alpha),\alpha}[\Vup(t+h,X^{t,x;u^S(\alpha),\alpha}_{t+h})] +\epsilon/2.
\end{align*}
Finally, by \eqref{ekv:G-stability} there is a $k'\geq k$ such that
\begin{align*}
  &\esssup_{u^S\in\mcU^{S,k'}_{t,t+h}}\essinf_{\alpha\in\mcA_{t,t+h}}G_{t,t+h}^{t,x;u^S(\alpha),\alpha}[ \Vup(t+h,X^{t,x;u^S(\alpha),\alpha}_{t+h})]
  \\
  &\leq \esssup_{u^S\in\mcU^{S,k'}_{t,t+h}}\essinf_{\alpha\in\mcA_{t,t+h}}G_{t,t+h}^{t,x;u^S(\alpha),\alpha}[ \Vup^{k'}(t+h,X^{t,x;u^S(\alpha),\alpha}_{t+h})]+\epsilon/2
\end{align*}
and we conclude that $\Vup$ satisfies \eqref{ekv:dynp-U}.\qed\\

%\begin{rem}\label{rem:dpp-at-stopping}
%Using the continuity of $\Vlow$ and $\Vup$ we see that the dynamic programming principles \eqref{ekv:dynp-W} and \eqref{ekv:dynp-U} still hold true if we replace $t+h$ with an arbitrary stopping time $\tau\in\mcT_t$.
%\end{rem}

\section{The value functions as viscosity solutions to the HJBI-QVI\label{sec:hjbi-qvi}}
Our main motivation for deriving the dynamic programming relations in the previous section is that we wish to use them to prove that the upper and lower value functions are solutions, in viscosity sense, to the Hamilton-Jacobi-Bellman-Isaacs quasi-variational inequality \eqref{ekv:var-ineq}.

Whenever $\Vlow(t,x)>\mcM \Vlow(t,x)$ (resp. $\Vup(t,x)>\mcM \Vup(t,x)$) a simple application of the dynamic programming principle stipulates that it is suboptimal for the impulse controller to intervene on the system at time $t$. One main ingredient when proving that $\Vlow$ (resp. $\Vup$) is a viscosity solution to \eqref{ekv:var-ineq} is showing that if $\Vlow(t,x)>\mcM \Vlow(t,x)$ (resp. $\Vup(t,x)>\mcM \Vup(t,x)$) then, on sufficiently small time intervals, we may (to a sufficient accuracy) assume that the impulse controller does not intervene on the system. As the probability that the state, when starting in $x$ at time $t$, leaves any ball with a finite radius containing $x$ on a non-empty interval $[t,t+h)$ is positive, this results requires a slightly intricate analysis compared to the deterministic setting (something that was pointed out already in \cite{TangHouSWgame}). In the following sequence of lemmas we extend the results from \cite{TangHouSWgame} to the case when the cost functional is defined in terms of the solution to a BSDE.

The first lemma is given without proof as it follows immediately from the definitions:
\begin{lem}\label{lem:mcM-monotone}
Let $u,v:[0,T]\times \R^n\to\R$ be locally bounded functions. $\mcM$ is monotone (if $u\leq v$ pointwise, then $\mcM u\leq \mcM v$). Moreover, $\mcM(u_*)$ (resp. $\mcM(u^*)$) is l.s.c.~(resp. u.s.c.).% and $(\mcM u)_*\leq
\end{lem}

In addition, rather than relying on the standard DPP from the previous section, formulated at deterministic times, we need the following ``weak'' DPP:
\begin{lem}\label{lem:alter-dpp}
Assume that $(t,x)\in[0,T]\times\R^n$ and $h\in[0,T-t]$, then for any $\alpha\in\mcA_{t,t+h}$ we have
\begin{align}
\Vlow(t,x)\leq\esssup_{\tau\in\mcT_t}G_{t,\tau\wedge t+h}^{t,x;\emptyset,\alpha}[\ett_{[\tau\leq t+h]}\mcM \Vlow(\tau,X^{t,x;\emptyset,\alpha}_{\tau})+\ett_{[\tau>t+h]}\Vlow(t+h,X^{t,x;\emptyset,\alpha}_{t+h})]\label{ekv:dynp2-W}
\end{align}
and
\begin{align}
\Vup(t,x)\leq\esssup_{\tau\in\mcT_t}G_{t,\tau\wedge t+h}^{t,x;\emptyset,\alpha}[\ett_{[\tau\leq t+h]}\mcM \Vup(\tau,X^{t,x;\emptyset,\alpha}_{\tau})+\ett_{[\tau>t+h]}\Vup(t+h,X^{t,x;\emptyset,\alpha}_{t+h})].\label{ekv:dynp2-U}
\end{align}
\end{lem}

\noindent \emph{Proof.} For $\epsilon>0$ we let $\alpha^S(u):=\alpha\oplus_{\tau_1}\alpha^{S,\epsilon}(u)$ for all $u\in\mcU_{t,t+h}$, where $\alpha^{S,\epsilon}\in\mcA_{t,t+h}$ is such that\footnote{We can repeat the approximation routine from Lemma~\ref{lem:W-lessthan-Wh} to show that such a strategy exists.}
\begin{align}%\nonumber
  &\ett_{[\tau_1\leq t+h]}\Vlow(\tau_1,X^{t,x;(\tau_1,\beta_1),\alpha}_{\tau_1})\geq \ett_{[\tau_1\leq t+h]}G_{\tau_1, t+h}^{\tau_1,X^{t,x;(\tau_1,\beta),\alpha}_{\tau_1};[u]_{2:},\alpha^{S,\epsilon}(u)}[\Vlow(t+h,X^{t,x;u,\alpha^{S}(u)}_{t+h})]-\epsilon,\label{ekv:weird-alpha}
\end{align}
$\Prob$-a.s., for all $u\in\mcU_{t,t+h}$, where $[u]_{2:}:=(\tau_i,\beta_i)_{2\leq i\leq N}$. Then $\alpha^S\in\mcA_{t,t+h}$ and there is a $u^\epsilon:=(\tau_i^\epsilon,\beta_i^\epsilon)_{1\leq i\leq N^\epsilon}\in\mcU_{t,t+h}$ such that
\begin{align*}
  \Vlow(t,x)\leq G_{t,t+h}^{t,x;u^\epsilon,\alpha^S(u^\epsilon)}[\Vlow(t+h,X^{t,x;u^\epsilon,\alpha^S(u^\epsilon)}_{t+h})]+\epsilon.
\end{align*}
On the other hand, the semi-group property of $G$ along with \eqref{ekv:weird-alpha} and comparison gives that
\begin{align*}
  &G_{t,t+h}^{t,x;u^\epsilon,\alpha^S(u^\epsilon)}[\Vlow(t+h,X^{t,x;u^\epsilon,\alpha^S(u^\epsilon)}_{t+h})]
  \\
  &= G_{t,\tau^\epsilon_1\wedge t+h}^{t,x;(\tau^\epsilon_1,\beta^\epsilon_1),\alpha^S(u^\epsilon)}[G_{\tau^\epsilon_1\wedge t+h,t+h}^{\tau^\epsilon_1,X^{t,x;(\tau^\epsilon_1,\beta^\epsilon_1),\alpha}_{\tau^\epsilon_1};[u^\epsilon]_{2:},\alpha^S(u^\epsilon)}[\Vlow(t+ h,X^{t,x;u^\epsilon,\alpha^S(u^\epsilon)}_{t+h})]]%-\ett_{[\tau^\epsilon_1\leq t+h]}\ell(\tau^\epsilon_1,X^{t,x;\emptyset,\alpha}_{\tau^\epsilon_1},\beta^\epsilon_1)
  \\
  &\leq G_{t,\tau^\epsilon_1\wedge t+h}^{t,x;\emptyset,\alpha}[\Vlow(\tau^\epsilon_1\wedge t+h ,X^{t,x;(\tau^\epsilon_1,\beta^\epsilon_1),\alpha}_{\tau^\epsilon_1\wedge t+h})-\ett_{[\tau^\epsilon_1\leq t+h]}\ell(\tau^\epsilon_1,X^{t,x;\emptyset,\alpha}_{\tau^\epsilon_1},\beta^\epsilon_1)+\epsilon]
  \\
  &\leq G_{t,\tau_1^\epsilon\wedge t+h}^{t,x;\emptyset,\alpha}[\ett_{[\tau_1^\epsilon\leq t+h]}\mcM \Vlow(\tau^\epsilon_1,X^{t,x;\emptyset,\alpha}_{\tau^\epsilon_1})+\ett_{[\tau_1^\epsilon> t+h]} \Vlow(\tau^\epsilon_1,X^{t,x;\emptyset,\alpha}_{\tau^\epsilon_1})]+C\epsilon.
\end{align*}
Since $\epsilon>0$ was arbitrary the first inequality follows by taking the essential supremum over all $\tau_1^\epsilon\in\mcT$.

Concerning the second inequality we have that for each $\epsilon>0$, there is a $u^{S}_\epsilon\in\mcU^S_{t,t+h}$ such that
\begin{align*}
  \Vup(t,x)\leq G_{t,t+h}^{t,x;u^S_\epsilon(\tilde \alpha),\tilde \alpha}[\Vup(t+h,X^{t,x;u^S_\epsilon(\tilde \alpha),\tilde \alpha}_{t+h})]+\epsilon
\end{align*}
for all $\tilde\alpha\in\mcA_{t,t+h}$. With $(\tau^\epsilon_1,\beta^\epsilon_1)=[u^S_\epsilon(\alpha)]_1$ (assuming that $\tau^\epsilon_1=\infty$ when $u^S_\epsilon(\alpha)$ does not contain interventions) we let $\alpha^\epsilon\in\mcA_{\tau^\epsilon_1,t+h}$ be such that
\begin{align*}
  \ett_{[\tau_1^\epsilon\leq t+h]}\Vup(\tau^\epsilon_1,X^{t,x;(\tau^\epsilon_1,\beta^\epsilon_1),\alpha}_{\tau^\epsilon_1})\geq \ett_{[\tau_1^\epsilon\leq t+h]}G_{\tau^\epsilon_1, t+h}^{\tau^\epsilon_1,X^{t,x;(\tau^\epsilon_1,\beta^\epsilon_1),\alpha}_{\tau^\epsilon_1}; [u^S_\epsilon(\alpha^\epsilon)]_{2:},\alpha^\epsilon}[ \Vup(t+h,X^{t,x;u^S_\epsilon(\alpha^\epsilon),\alpha}_{t+h})]-\epsilon.%\label{ekv:weird-u}
\end{align*}
Applying the continuous control $\alpha\oplus_{\tau^\epsilon_1}\alpha^\epsilon$ (see Remark~\ref{rem:@tph} concerning the value at the point of concatenation) and using the semi-group property of $G$ (as above) now leads to the second inequality.\qed\\

\begin{lem}\label{lem:do-nothing}
Let $(t,x)\in [t,T)\times \R^n$ be such that $\Vlow(t,x)>\mcM \Vlow(t,x)$ then there is a $C>0$ and an $h'\in (0,T-t]$ such that
\begin{align*}
  \Vlow(t,x)\leq \essinf_{\alpha\in\mcA_{t,t+h}}G_{t,t+h}^{t,x;\emptyset,\alpha}[\Vlow(t+h,X^{t,x;\emptyset,\alpha}_{t+h})]+ Ch^{3/2}
\end{align*}
for all $h\in[0,h']$.
\end{lem}

\noindent \emph{Proof.} Since $\Vlow$ is continuous, Lemma~\ref{lem:mcM-monotone} implies that so is $\mcM\Vlow$. There is thus a $h''>0$ and an $\epsilon>0$ such that
\begin{align*}
  \inf_{(t',x')\in [t,t+h'']\times \bar B_\epsilon(x)}\Vlow(t',x')\geq \sup_{(t',x')\in [t,t+h'']\times \bar B_\epsilon(x)}\mcM \Vlow(t',x')+\epsilon,
\end{align*}
with $\bar B_\epsilon(x):=\{x'\in\R^n:|x'-x|\leq\epsilon\}$.

For each $\alpha\in\mcA_{t,t+h}$ we have, with $X:=X^{t,x;\emptyset,\alpha}$, for $p\geq 2$ that
\begin{align}\nonumber
\E\Big[\sup_{s\in [t,t+h]}|X_s-x|^p\Big|\mcF_t\Big]&\leq C\E\Big[(\int_t^{t+h}|a(s,X^{t,x}_s)|ds)^p+(\int_t^{t+h}|\sigma(s,X_s)|^2ds)^{p/2}\Big|\mcF_t\Big]
\\
&\leq Ch^{p/2}(1+\E\Big[\sup_{s\in [t,t+h]}|X_s|^p\Big|\mcF_t\Big])\nonumber
\\
&\leq Ch^{p/2}(1+|x|^p),\label{ekv:divergence-bound}
\end{align}
$\Prob$-a.s. We introduce the stopping time
\begin{align*}
  \eta:=\inf\big\{s\geq t: X_s\notin B_\epsilon(x)\big\}
\end{align*}
(with $\inf\emptyset=\infty$) and get that
\begin{align*}
\E\big[\ett_{[\eta\leq t+h]}\big|\mcF_t \big]\epsilon^p\leq\E\Big[\sup_{s\in [t,t+h]}|X_s-x|^p\Big|\mcF_t\Big]&\leq Ch^{p/2}(1+|x|^p),
\end{align*}
$\Prob$-a.s. Choosing $p=6$ gives
\begin{align*}
\E\big[\ett_{[\eta\leq t+h]}\big|\mcF_t \big]\leq \epsilon^{-6}Ch^{3}(1+|x|^6),
\end{align*}
$\Prob$-a.s. Using this inequality we will show that there is a $C>0$ such that for some $h'\in(0,h'']$ and all $h\in [0,h']$ we have
\begin{align*}
\esssup_{\tau\in\mcT_t}G_{t,\tau\wedge t+h}^{t,x;\emptyset,\alpha}[\ett_{[\tau\leq t+h]}\mcM \Vlow(\tau,X_{\tau})+\ett_{[\tau>t+h]}\Vlow(t+h,X_{t+h})]\leq G_{t,t+h}^{t,x;\emptyset,\alpha}[\Vlow(t+h,X_{t+h})]+ Ch^{3/2}
\end{align*}
for all $\alpha\in\mcA_{t,t+h}$ from which the result of this lemma follows by Lemma~\ref{lem:alter-dpp}. For any $\tau\in\mcT_t$, let $(\mcY^1,\mcZ^1)$ be the unique solution to
\begin{align*}
\mcY^1_s&=\ett_{[\tau\leq t+h]}\mcM \Vlow(\tau,X_{\tau}) + \ett_{[\tau>t+h]}\Vlow(t+h,X_{t+h})
\\
&\quad+\int_s^{\tau\wedge t+h} f(r,X_r,\mcY^1_r,\mcZ^1_r,\alpha_r)dr-\int_s^{\tau\wedge t+h} \mcZ^1_rdW_r.
\end{align*}
with $X:=X^{t,x;\emptyset,\alpha}$ and let $(\mcY^2,\mcZ^2)$ solve
\begin{align*}
\mcY^2_s&=\Vlow(t+h,X_{t+h})+\int_s^{t+h} f(r,X_r,\mcY^2_r,\mcZ^2_r,\alpha_r)dr-\int_s^{t+h} \mcZ^2_rdW_r.
\end{align*}
Then, with
\begin{align}\label{ekv:zeta_1-in-DN}
 \zeta_1(s):=\frac{f(s,X_s,\mcY^2_s,\mcZ^2_s,\alpha_s) -f(s,X_s,\ett_{[s\leq\tau]}\mcY^1_s,\mcZ^2_s,\alpha_s)} {\mcY^2_s-\mcY^1_s} \ett_{[\mcY^2_s\neq\ett_{[s\leq\tau]}\mcY^1_s]}
\end{align}
and
\begin{align}\label{ekv:zeta_2-in-DN}
 \zeta_2(s):=\frac{f(s,X_s,\ett_{[s\leq\tau]}\mcY^1_s,\mcZ^2_s,\alpha_s) -f(s,X_s,\ett_{[s\leq\tau]}\mcY^1_s,\ett_{[s\leq\tau]}\mcZ^1_s,\alpha_s)} {|\mcZ^2_s-\ett_{[s\leq\tau]}\mcZ^1_s|^2} (\mcZ^2_s-\ett_{[s\leq\tau]}\mcZ^1_s)^\top
\end{align}
we have by the Lipschitz continuity of $f$ that $|\zeta_1(s)|\vee|\zeta_2(s)|\leq k_f$. Then, with
\begin{align*}
  M_{r,s}:=e^{\int_r^{s}(\zeta_1(v)-\frac{1}{2}|\zeta_2(v)|^2)dv+\frac{1}{2}\int_r^{s}\zeta_2(v)dW_v}=:e^{\int_r^{s}\zeta_1(v)dv}\tilde M_{r,s},
\end{align*}
we have
\begin{align*}
\mcY^1_t-\mcY^2_t&=\E\Big[\ett_{[\tau< t+h]}\Big\{M_{t,\tau}(\mcM \Vlow(\tau,X_{\tau})-M_{\tau,t+h}\Vlow(t+h,X_{t+h}))
-\int_\tau^{t+h}M_{t,s}f(s,X_s,0,0)ds\Big\}\Big|\mcF_t\Big]
\\
&=\Lambda_1(h)+\Lambda_2(h),
\end{align*}
where
\begin{align*}
\Lambda_1(h):\!&\!=\E\Big[\ett_{[\eta<t+h]}\ett_{[\tau< t+h]}\Big\{M_{t,\tau}(\mcM \Vlow(\tau,X_{\tau})-M_{\tau,t+h}\Vlow(t+h,X_{t+h}))
\\
&\quad-\int_\tau^{t+h}M_{t,s}f(s,X_s,0,0)ds\Big\}\Big|\mcF_t\Big]
\\
&\leq C\E\big[\ett_{[\eta<t+h]}\big]^{1/2}\E\Big[|M_{t,\tau}\mcM \Vlow(\tau,X^{t,x;\emptyset,\alpha}_{\tau})|^2+|M_{t,t+h}\Vlow(t+h,X_{t+h})|^2
\\
&\quad+\int_\tau^{t+h}|M_{t,s}f(s,X_s,0,0)|^2ds\Big|\mcF_t\Big]^{1/2}
\\
&\leq C(1+|x|^\rho)h^{3/2}
\end{align*}
and
\begin{align*}
\Lambda_2(h):\!&\!=\E\Big[\ett_{[\eta\geq t+h]}\ett_{[\tau< t+h]}\Big\{M_{t,\tau}(\mcM \Vlow(\tau,X_{\tau})-M_{\tau,t+h}\Vlow(t+h,X_{t+h}))
\\
&\quad-\int_\tau^{t+h}M_{t,s}f(s,X_s,0,0)ds\Big\}\Big|\mcF_t\Big]
\\
&\leq \E\Big[\ett_{[\eta\geq t+h]}\ett_{[\tau< t+h]}\Big\{M_{t,\tau}(\mcM \Vlow(\tau,X_{\tau})-\Vlow(t+h,X_{t+h})+(1-M_{\tau,t+h})\Vlow(t+h,X_{t+h}))
\\
&\quad+C(1+|x|^\rho)\int_{\tau}^{t+h}M_{t,s}ds\Big\}\Big|\mcF_t\Big]
\\
&\leq \E\Big[\ett_{[\eta\geq t+h]}\ett_{[\tau< t+h]}\Big\{-\epsilon M_{t,\tau}+(M_{t,t+h}-M_{t,\tau}+\int_{\tau}^{t+h}M_{t,s}ds)C(1+|x|^\rho)\Big\}\Big|\mcF_t\Big],
\end{align*}
where we have used the polynomial growth of $\Vlow$ and $\mcM \Vlow$ together with the fact that $\sup_{s\in[t,t+h]}|X_s|\leq |x|+\epsilon$ on $\{\omega:\eta\geq t+h\}$. We can now get rid of $\ett_{[\eta\geq t+h]}$ and use the martingale property of $\tilde M$ to find that
\begin{align*}
\Lambda_2(h)&\leq \E\Big[\ett_{[\tau< t+h]}\Big\{-\epsilon M_{t,\tau}+(M_{t,t+h}-M_{t,\tau}+\int_{\tau}^{t+h}M_{t,s}ds)C(1+|x|^\rho)\Big\}\Big|\mcF_t\Big]+C(1+|x|^\rho)h^{3/2}
\\
&\leq\E\Big[\ett_{[\tau< t+h]}e^{\int_t^{\tau}\zeta_1(r)dr}\tilde M_{t,t+h}\Big\{-\epsilon +C(1+|x|^\rho)(e^{k_f h}-1+he^{k_f h})
\Big\}\Big|\mcF_t\Big]+C(1+|x|^\rho)h^{3/2}
\\
&\leq C(1+|x|^\rho)h^{3/2}
\end{align*}
whenever $h>0$ is small enough that $-\epsilon +C(1+|x|^\rho)(e^{k_f h}-1+he^{k_f h})\leq 0$. Combined, this gives that there is a $h'\in [0,h'']$ such that whenever $h\in [0,h']$ we have $\mcY^1_t-\mcY^2_t\leq C(1+|x|^\rho)h^{3/2}$. Since $\tau$ and $\alpha$ were arbitrary the assertion follows.\qed\\

\begin{lem}\label{lem:do-nothing-U}
Let $(t,x)\in [t,T)\times \R^n$ be such that $\Vup(t,x)>\mcM \Vup(t,x)$ then there is a $C>0$ and an $h'\in (0,T-t]$ such that
\begin{align*}
  \Vup(t,x)\leq \essinf_{\alpha\in\mcA_{t,t+h}}G_{t,t+h}^{t,x;\emptyset,\alpha}[\Vup(t+h,X^{t,x;\emptyset,\alpha}_{t+h})]+ Ch^{3/2}
\end{align*}
for all $h\in[0,h')$.
\end{lem}

\noindent\emph{Proof.} As in the proof of the above lemma, there is a $h''>0$ and an $\epsilon>0$ such that
\begin{align*}
  \inf_{(t',x')\in [t,t+h'']\times \bar B_\epsilon(x)}\Vup(t',x')\geq \sup_{(t',x')\in [t,t+h'']\times \bar B_\epsilon(x)}\mcM \Vup(t',x')+\epsilon,
\end{align*}
We can thus repeat the steps in the previous lemma to conclude that there is a $C>0$ such that
\begin{align*}
\esssup_{\tau\in\mcT_t}G_{t,\tau\wedge t+h}^{t,x;\emptyset,\alpha} [\ett_{[\tau\leq t+h]}\mcM \Vup(\tau,X_{\tau})+\ett_{[\tau>t+h]}\Vup(t+h,X_{t+h})]\leq G_{t,t+h}^{t,x;\emptyset,\alpha}[\Vup(t+h,X_{t+h})]+ Ch^{3/2}
\end{align*}
for all $\alpha\in\mcA_{t,t+h}$ and $h\in [0,h']$ for some $h'\in (0,h'']$. The lemma then follows by applying the second inequality in Lemma~\ref{lem:alter-dpp}.\qed\\

We now fix $(t,x)\in[0,T)\times \R^n$, $h\in (0,T-t]$ and $\varphi\in C^{3}_{l,b}$. Following the standard procedure to go from a DPP to a quasi-variational inequality when dealing with a controlled FBSDE (see \eg \cite{Li14}) we introduce the BSDEs
\begin{align*}
  Y^{1,\alpha}_{s}=\int_s^{t+h}F(r,X^{t,x;\emptyset,\alpha}_r,Y^{1,\alpha}_r,Z^{1,\alpha}_r,\alpha_r)dr-\int_s^{t+h}Z^{1,\alpha}_rdW_r
\end{align*}
and
\begin{align}\label{ekv:bsde-Y2}
  Y^{2,\alpha}_{s}=\int_s^{t+h}F(r,x,Y^{2,\alpha}_r,Z^{2,\alpha}_r,\alpha_r)dr-\int_s^{t+h}Z^{2,\alpha}_rdW_r,
\end{align}
with
\begin{align*}
  F(s,x,y,z,\alpha)&:=\frac{\partial}{\partial s}\varphi(s,x)+\frac{1}{2}\trace\{\sigma\sigma^\top(s,x,\alpha)D^2_x\varphi(s,x)\}+(D_x\varphi(s,x))b(s,x,\alpha)
  \\
  &\quad +f(s,x,\varphi(s,x)+y,(D_x\varphi(s,x))\sigma(s,x,\alpha)+z,\alpha).
\end{align*}

\begin{rem}
It is easy to check that the driver $F$ satisfies Assumption~\ref{ass:on-coeff}.\ref{ass:on-coeff-f} from which we conclude that the above BSDEs both admit unique solutions.
\end{rem}

In particular, we note that $u$ is a viscosity supersolution (subsolution) of \eqref{ekv:var-ineq} if $u(T,x)\geq(\leq)\psi(T,x)$, $u(t,x)\geq \mcM u(t,x)$ and $\inf_{\alpha\in A}F(t,x,0,0,\alpha)\leq 0$ ($\geq 0$) on $\mcD_C(u):=\{(t,x):u(t,x)> \mcM u(t,x)\}$ whenever $\varphi\in C^3_{l,b}$ is such that $u(t,x)=\varphi(t,x)$ and $u(t,x)-\varphi(t,x)$ attains a local minimum (maximum) at $(t,x)$.%\footnote{We may w.l.o.g.~assume that $\varphi$ in Definition~\ref{def:visc-sol-dom} belongs to $C^3_b$.}

Note that the only reason that \eqref{ekv:bsde-Y2} is stochastic comes from the fact that $\alpha$ is a stochastic control. In regard to Hamiltonian minimization it seems natural to introduce the following ordinary differential equation (ODE)
\begin{align*}
  Y^{0}_{s}=\int_s^{t+h}\inf_{\alpha\in A}F(s,x,Y^{0}_s,0,\alpha)ds.
\end{align*}

We have the following auxiliary lemma, that summarize the results in Lemma 5.1 and Lemma 5.3 of \cite{Buckdahn08}.

\begin{lem}\label{lem:Y1-to-G}
For every $\alpha\in\mcA_{t,t+h}$ and $s\in[t,t+h]$ we have
\begin{align}\label{ekv:Y1-to-G}
  Y^{1,\alpha}_s=G^{t,x;\emptyset,\alpha}_{s,t+h}[\varphi(t+h,X^{t,x;\emptyset,\alpha}_{t+h})]-\varphi(s,X^{t,x;\emptyset,\alpha}_{s}),\qquad \Prob-{\rm a.s.}
\end{align}
Also, we have that
\begin{align}\label{ekv:Y0-rep}
  Y^{0}_t=\essinf_{\alpha\in\mcA_{t,t+h}}Y^{2,\alpha}_t,\qquad \Prob-{\rm a.s.}
\end{align}
\end{lem}

\noindent \emph{Proof.} The first property follows from the definition of $G$ and Ito's formula applied to $\varphi(s,X^{t,x;\emptyset,\alpha}_{s})$. The second result is immediate from the comparison principle of BSDEs.\qed\\

We now give a sequence of lemmata that will help us show that $\Vlow$ is a viscosity solution to \eqref{ekv:var-ineq}.

\begin{lem}\label{lem:Y1-Y2-diff}
We have
\begin{align*}
  |Y^{1,\alpha}_t-Y^{2,\alpha}_t|\leq Ch^{3/2},\qquad \Prob-{\rm a.s.}
\end{align*}
\end{lem}

\noindent \emph{Proof.}  Note that
\begin{align*}
|Y^{1,\alpha}_t-Y^{2,\alpha}_t|&\leq\E\Big[\int_t^{t+h}| F(s,X^{t,x;\emptyset,\alpha}_s,Y^{1,\alpha}_s,Z^{1,\alpha}_s,\alpha_s)-F(s,x,Y^{2,\alpha}_s,Z^{2,\alpha}_s,\alpha_s)|ds \Big|\mcF_t\Big]
\\
&\leq C\E\Big[\int_t^{t+h}((1+|x|^\rho+|X^{t,x;\emptyset,\alpha}_s|^\rho)|X^{t,x;\emptyset,\alpha}_s-x|+|Y^{1,\alpha}_s-Y^{2,\alpha}_s|+|Z^{1,\alpha}_s -Z^{2,\alpha}_s|)ds \Big|\mcF_t\Big].
\end{align*}
Concerning the first term on the right-hand side we have
\begin{align*}
\E\Big[\int_t^{t+h}(1+|x|^\rho+|X^{t,x;\emptyset,\alpha}_s|^\rho)|X^{t,x;\emptyset,\alpha}_s-x|ds \Big|\mcF_t\Big]&\leq C\E\Big[\sup_{s\in[t,t+h]}|X^{t,x;\emptyset,\alpha}_s-x|^{2}\Big|\mcF_t\Big]^{1/2}h
\\
&\leq Ch^{3/2}
\end{align*}%+\E\Big[\sup_{s\in[t,t+h]}|X^{t,x;\emptyset,\alpha}_s|^{2\rho}\Big|\mcF_t\Big]^{1/2})
For the remaining terms,
\begin{align*}
&\E\Big[\int_t^{t+h}(|Y^{1,\alpha}_s-Y^{2,\alpha}_s|+|Z^{1,\alpha}_s-Z^{2,\alpha}_s|)ds \Big|\mcF_t\Big]
\\
&\leq \sqrt 2\E\Big[\int_t^{t+h}(|Y^{1,\alpha}_s-Y^{2,\alpha}_s|^2+|Z^{1,\alpha}_s-Z^{2,\alpha}_s|^2)ds \Big|\mcF_t\Big]^{1/2}h^{1/2}
\end{align*}
and classically we have
\begin{align*}
&\E\Big[\int_t^{t+h}|Y^{1,\alpha}_s-Y^{2,\alpha}_s|^2+|Z^{1,\alpha}_s-Z^{2,\alpha}_s|^2)ds \Big|\mcF_t\Big]
\\
&\leq C\E\Big[\int_t^{t+h}|F(s,X^{t,x;\emptyset,\alpha}_s,Y^{1,\alpha}_s,Z^{1,\alpha}_s,\alpha_s)-F(s,x,Y^{1,\alpha}_s,Z^{1,\alpha}_s,\alpha_s)|^2ds \Big|\mcF_t\Big]
\\
&\leq C\E\Big[\int_t^{t+h}(1+|x|^\rho+|X^{t,x;\emptyset,\alpha}_s|^{2\rho})|X^{t,x;\emptyset,\alpha}_s-x|^2ds \Big|\mcF_t\Big]
\\
&\leq Ch^{2}.
\end{align*}
Combining the above estimates the desired results follows.\qed\\

\begin{lem}\label{lem:Y0-int}
There is a $C>0$ such that
\begin{align*}
  \int_t^{t+h}|Y^{0}_s|ds\leq Ch^{3/2}.
\end{align*}
for each $t\in[0,T]$ and $h\in [0,T-t]$.
\end{lem}

\noindent \emph{Proof.} Gr\"onwall's inequality gives that
\begin{align*}
\sup_{s\in [t,t+h]}|Y^{0,\alpha}_s|\leq Ch
\end{align*}
and we conclude that $\int_t^{t+h}|Y^{0,\alpha}_s|ds\leq h\sup_{s\in [t,t+h]}|Y^{0,\alpha}_s|\leq Ch^2$.
\qed\\
%Moreover, by definition
%\begin{align*}
%|Y^{2,\alpha}_s|&\leq \E\Big[\int_s^{t+h}|F(r,x,Y^{2,\alpha}_r,Z^{2,\alpha}_r,\alpha_r)|dr\Big|\mcF_s\Big]
%\\
%&\leq C\E\Big[\int_s^{t+h}(1+|x|^{\rho}+|Y^{2,\alpha}_r|+|Z^{2,\alpha}_r|)dr\Big|\mcF_s\Big]
%\\
%&\leq C(h+h^{1/2}\E\Big[\int_s^{t+h}|Z^{2,\alpha}_r|^2dr\Big|\mcF_s\Big]^{1/2})
%\\
%&\leq Ch.
%\end{align*}
%Now, since $Y^{2,\alpha}_s=\int_s^{t+h}F(r,x,Y^{2,\alpha}_r,Z^{2,\alpha}_r)dr-\int_s^{t+h}Z^{2,\alpha}_rdW_r$ we find that
%\begin{align*}
%\E\Big[\int_t^{t+h}|Z^{2,\alpha}_s|^2ds\Big|\mcF_t\Big]&\leq C\E\Big[|Y^{2,\alpha}_t|^2+h\int_t^{t+h}|F(r,x,Y^{2,\alpha}_r,Z^{2,\alpha}_r,\alpha_r)|^2dr\Big|\mcF_t\Big]
%\\
%&\leq C\E\Big[|Y^{2,\alpha}_t|^2+h^2+h\int_t^{t+h}(|Y^{2,\alpha}_r|^2+|Z^{2,\alpha}_r|^2)dr\Big|\mcF_t\Big],
%\end{align*}
%whence $\E\Big[\int_t^{t+h}|Z^{2,\alpha}_s|^2ds\Big|\mcF_t\Big]\leq Ch^2$ for $h<h_0$ for some $h_0>0$. By choosing $C>0$ sufficiently large and using \eqref{ekv:YZmom1} this extends to all $h\in [0,T-t]$. We thus conclude that
%\begin{align*}
%\E\Big[\int_t^{t+h}|Z^{2,\alpha}_s|ds\Big|\mcF_t\Big]\leq h^{1/2}\E\Big[\int_t^{t+h}|Z^{2,\alpha}_s|^2ds\Big|\mcF_t\Big]^{1/2}\leq Ch^{3/2}
%\end{align*}
%and the results follows.\qed\\

\begin{lem}\label{lem:local-minimum}
Assume that $\varphi\in \Pi_{pg}$ is such that $\varphi-\Vlow$ has a local maximum at $(t,x)$ where $\varphi(t,x)=\Vlow(t,x)$. Then, there are constants $C,h'>0$ such that
\begin{align*}
G_{t,t+h}^{t,x;\emptyset,\alpha}[\Vlow(t+h,X^{t,x;\emptyset,\alpha}_{t+h})]\geq G_{t,t+h}^{t,x;\emptyset,\alpha}[\varphi(t+h,X^{t,x;\emptyset,\alpha}_{t+h})]-Ch^{3/2}
\end{align*}
for all $h\in [0,(T-t)\wedge h']$ and $\alpha\in\mcA_{t,t+h}$.
\end{lem}

\noindent \emph{Proof.} Since $\varphi-\Vlow$ has a local maximum at $(t,x)$ there are constants $\epsilon,h'>0$ and a $h'>0$ such that $\Vlow(t',x')\geq \varphi(t',x')$ for all $(t',x')\in [t,t+h'\wedge T]\times \bar B_\epsilon(x)$. Now, let
\begin{align*}
  \eta:=\inf\{s\geq t: X^{t,x;\emptyset,\alpha}\notin B_\epsilon(x)\}
\end{align*}
and note from the proof of Lemma~\ref{lem:do-nothing} that $\E[\ett_{[\eta\leq t+h]}|\mcF_t]\leq Ch^3$, $\Prob$-a.s. Assume that $h\in [0,T-t]$ and let $(\mcY^1,\mcZ^1)$ be the unique solution to
\begin{align*}
\mcY^1_s&=\Vlow(t+h,X^{t,x;\emptyset,\alpha}_{t+h})+\int_s^{t+h} f(r,X^{t,x;\emptyset,\alpha}_r,\mcY^1_r,\mcZ^1_r,\alpha_r)dr-\int_s^{t+h} \mcZ^1_rdW_r.
\end{align*}
and assume that $(\mcY^2,\mcZ^2)$ satisfies
\begin{align*}
\mcY^2_s&=\varphi(t+h,X^{t,x;\emptyset,\alpha}_{t+h})+\int_s^{t+h} f(r,X^{t,x;\emptyset,\alpha}_r,\mcY^2_r,\mcZ^2_r,\alpha_r)dr-\int_s^{t+h} \mcZ^2_rdW_r.
\end{align*}
Then, with
\begin{align*}
  M_s:=e^{\int_t^{s}(\zeta_1(r)-\frac{1}{2}\zeta_2^2(r))dr+\frac{1}{2}\int_t^{s}\zeta_2(r)dW_r},
\end{align*}
where $\zeta_1$ and $\zeta_2$ are given by \eqref{ekv:zeta_1-in-DN}-\eqref{ekv:zeta_2-in-DN}. By comparison we have
\begin{align*}
\mcY^2_t-\mcY^1_t&\leq \E\big[\ett_{[\eta\leq t+h]}M_{t+h}(\varphi(t+h,X^{t,x;\emptyset,\alpha}_{t+h})-\Vlow(t+h,X^{t,x;\emptyset,\alpha}_{t+h}))\big|\mcF_t\big]
\\
&\leq \sqrt{2}\E\big[\ett_{[\eta\leq t+h]}\big|\mcF_t\big]^{1/2}\E\big[M_{t+h}^2(|\varphi(t+h,X^{t,x;\emptyset,\alpha}_{t+h})|^2 +|\Vlow(t+h,X^{t,x;\emptyset,\alpha}_{t+h})|^2)\big|\mcF_t\big]^{1/2}
\\
&\leq Ch^{3/2}
\end{align*}
and the result follows.\qed\\

\begin{thm}
$\Vlow$ is a viscosity solution to \eqref{ekv:var-ineq}.
\end{thm}

\noindent \emph{Proof.} To begin with we clearly have that $\Vlow(T,x)=\psi(x)$ for all $x\in\R^n$ (see Remark~\ref{rem:@end}). We first show that $\Vlow$ is a viscosity supersolution. For this, we fix $(t,x)\in [0,T]\times \R^n$ and assume that $\varphi$ is such that $\varphi-\Vlow$ has a local maximum at $(t,x)$, where $\varphi(t,x)=\Vlow(t,x)$.

If $(t,x)\in \mcD_C(\Vlow)$ we have by the DPP that %:=\{(t',x')\in [0,T]\times \R^n:\Vlow(t',x')>\mcM \Vlow(t',x')\}
\begin{align*}
\varphi(t,x)=\Vlow(t,x)&=\essinf_{\alpha^S\in\mcA^S_{t,t+h}}\esssup_{u\in\mcU^k_{t,t+h}}G_{t,t+h}^{t,x;u,\alpha^S(u)}[\Vlow(t + h,X^{t,x;u,\alpha^S(u)}_{t+h})]
\\
&\geq \essinf_{\alpha\in\mcA_{t,t+h}}G_{t,t+h}^{t,x;\emptyset,\alpha}[\Vlow(t+h,X^{t,x;\emptyset,\alpha}_{t+h})]
\end{align*}
On the other hand by Lemma~\ref{lem:local-minimum} we have for $h>0$ sufficiently small that
\begin{align*}
G_{t,t+h}^{t,x;\emptyset,\alpha}[\varphi(t+h,X^{t,x;\emptyset,\alpha}_{t+h})]\leq G_{t,t+h}^{t,x;\emptyset,\alpha}[\Vlow(t+h,X^{t,x;\emptyset,\alpha}_{t+h})]+Ch^{3/2}.
\end{align*}
Now, \eqref{ekv:Y1-to-G} gives
\begin{align*}
Y^{1,\alpha}_t\leq G_{t,t+h}^{t,x;\emptyset,\alpha}[\Vlow(t+h,X^{t,x;\emptyset,\alpha}_{t+h})]-\varphi(t,x)+Ch^{3/2}.
\end{align*}
Combined this gives
\begin{align*}
\essinf_{\alpha\in\mcA_{t,t+h}}Y^{1,\alpha}_t\leq Ch^{3/2}.
\end{align*}
In particular, by Lemma~\ref{lem:Y1-Y2-diff} and \eqref{ekv:Y0-rep} this implies that
\begin{align*}
Y^0_t=\essinf_{\alpha\in\mcA_{t,t+h}}Y^{2,\alpha}_t\leq Ch^{3/2}.
\end{align*}
Hence, $\lim_{h\to 0} h^{-1}Y^0_t\leq 0$ and we conclude by Lemma~\ref{lem:Y0-int} that
\begin{align*}
  0&\geq\lim_{h\to 0} h^{-1}\int_{t}^{t+h}\inf_{\alpha\in A}F(s,x,Y^0_s,0,\alpha)ds
  \\
  &\geq \lim_{h\to 0} h^{-1}\int_{t}^{t+h}\inf_{\alpha\in A}(F(s,x,0,0,\alpha)-C|Y^0_s|)ds
  \\
  &=\lim_{h\to 0} h^{-1}\int_{t}^{t+h}\inf_{\alpha\in A}F(s,x,0,0,\alpha)ds
\end{align*}
and by continuity of $\inf_{\alpha\in A}F(\cdot,x,0,0,\alpha)$ it follows that
\begin{align*}
  \inf_{\alpha\in A} F(t,x,0,0,\alpha)\leq 0.
\end{align*}

Assume instead that $(t,x)\in \mcD_S(\Vlow):=([0,T]\times\R^n)\setminus \mcD_C(\Vlow)$, then $\Vlow(t,x)=\mcM \Vlow(t,x)$ and we conclude that $\Vlow$ is a viscosity supersolution.\\

We turn now to the subsolution property. We fix $(t,x)\in [0,T]\times \R^n$ and assume that $\varphi$ is such that $\varphi-\Vlow$ has a local minimum at $(t,x)$, where $\varphi(t,x)=\Vlow(t,x)$.
If $(t,x)\in \mcD_C(\Vlow)$ we have by the DPP and Lemma~\ref{lem:do-nothing} that, whenever $h>0$ is sufficiently small,
\begin{align*}
\varphi(t,x)=\Vlow(t,x)&=\essinf_{\alpha^S\in\mcA^S_{t,t+h}}\esssup_{u\in\mcU^k_{t,t+h}}G_{t,t+h}^{t,x;u,\alpha^S(u)}[\Vlow(t + h,X^{t,x;u,\alpha^S(u)}_{t+h})]
\\
&\leq \essinf_{\alpha\in\mcA_{t,t+h}}G_{t,t+h}^{t,x;\emptyset,\alpha}[\Vlow(t+h,X^{t,x;\emptyset,\alpha}_{t+h})] + Ch^{3/2}
\end{align*}
On the other hand repeating the argument in the proof of Lemma~\ref{lem:local-minimum} gives that
\begin{align*}
G_{t,t+h}^{t,x;\emptyset,\alpha}[\Vlow(t+h,X^{t,x;\emptyset,\alpha}_{t+h})]\leq G_{t,t+h}^{t,x;\emptyset,\alpha}[\varphi(t+h,X^{t,x;\emptyset,\alpha}_{t+h})]+ Ch^{3/2}
\end{align*}
and we get that
\begin{align*}
-Y^{1,\alpha}_t&=\varphi(t,x)-G_{t,t+h}^{t,x;\emptyset,\alpha^S(\emptyset)}[\varphi(t+h,X^{t,x;\emptyset,\alpha^S(\emptyset)}_{t+h})]
\\
&\leq Ch^{3/2},
\end{align*}
\ie $Y^{1,\alpha}_t\geq -Ch^{3/2}$. Now, repeating the above argument we find that
\begin{align*}
  \inf_{\alpha\in A} F(t,x,0,0,\alpha)\geq 0.
\end{align*}
Analogously we get when $(t,x)\in \mcD_S(\Vlow)$ then $\Vlow(t,x)= \mcM\varphi(t,x)$ and we conclude that $\Vlow$ is a viscosity subsolution.\qed\\

\begin{rem}
By the same argument while using Lemma~\ref{lem:do-nothing-U} instead of Lemma~\ref{lem:do-nothing} we conclude that $\Vup$ is a viscosity solution to \eqref{ekv:var-ineq}.
\end{rem}

\section{Uniqueness of viscosity solutions to the HJBI-QVI\label{sec:unique}}

To be able to conclude that the game has a value, \ie that $\Vlow\equiv\Vup$, we will now show that \eqref{ekv:var-ineq} has at most one solution in the viscosity sense in $\Pi_{pg}$. We let
\begin{align}
  \mcL^\alpha\varphi(t,x):=\sum_{j=1}^d a_j(t,x,\alpha)\frac{\partial}{\partial x_j}\varphi(t,x)+\frac{1}{2}\sum_{i,j=1}^d (\sigma\sigma^\top(t,x,\alpha))_{i,j}\frac{\partial^2}{\partial x_i\partial x_j}\varphi(t,x)
\end{align}
and have that
\begin{align*}
  H(t,x,v(t,x),Dv(t,x),D^2v(t,x),\alpha):=\mcL^\alpha v(t,x)+f(t,x,v(t,x),Dv(t,x)\cdot \sigma(t,x,\alpha),\alpha).
\end{align*}

We will need the following lemma:

\begin{lem}\label{lem:is-super}
Let $v$ be a supersolution to \eqref{ekv:var-ineq} satisfying
\begin{align*}
\forall (t,x)\in[0,T]\times\R^d,\quad |v(t,x)|\leq C(1+|x|^{2\gamma})
\end{align*}
for some $\gamma>0$. Then there is a $\lambda_0 > 0$ such that for any $\lambda>\lambda_0$ and $\theta> 0$, the function $v + \theta e^{-\lambda t}(1+((|x|-K_\Gamma)^+)^{2\gamma+2})$ is also a supersolution to \eqref{ekv:var-ineq}.
\end{lem}

\noindent \emph{Proof.} With $w:=v + \theta e^{-\lambda t}(1+((|x|-K_\Gamma)^+)^{2\gamma+2})$ we note that, since $v$ is a supersolution and $\theta e^{-\lambda T}(1+((|x|-K_\Gamma)^+)^{2\gamma+2})\geq 0$, we have $w(T,x)\geq v(T,x)\geq \psi(x)$ so that the terminal condition holds. Moreover, we have
\begin{align*}
&w(t,x)-\sup_{b\in U}\{w(t,\Gamma(t,x,b))-\ell(t,x,b)\}
\\
&=v(t,x) + \theta e^{-\lambda t}(1+((|x|-K_\Gamma)^+)^{2\gamma+2})
\\
&\quad-\sup_{b\in U}\{v(t,\Gamma(t,x,b)) + \theta e^{-\lambda t}(1+((|\Gamma(t,x,b)|-K_\Gamma)^+)^{2\gamma+2}-\ell(t,x,b))\}
\\
&\geq v(t,x)- \sup_{b\in U}\{v(t,\Gamma(t,x,b))-\ell(t,x,b)\}
\\
&\quad+\theta e^{-\lambda t}\{(1+((|x|-K_\Gamma)^+)^{2\gamma+2})- \sup_{b\in U}(1+((|\Gamma(t,x,b)|-K_\Gamma)^+)^{2\gamma+2})\}.
\end{align*}
Since $v$ is a supersolution, we have
\begin{align*}
  v(t,x)- \sup_{b\in U}\{v(t,\Gamma(t,x,b))-\ell(t,x,b)\}\geq 0
\end{align*}
Now, either $|x|\leq K_\Gamma$ in which case it follows by \eqref{ekv:imp-bound} that $|\Gamma(t,x,b)|\leq K_\Gamma$ or $|x|> K_\Gamma$ and \eqref{ekv:imp-bound} gives that $|\Gamma(t,x,b)|\leq |x|$. We conclude that
\begin{align*}
  w(t,x)- \sup_{b\in U}\{w(t,\Gamma(t,x,b))-\ell(t,x,b)\}\geq 0.
\end{align*}
Next, let $\varphi\in C^{1,2}([0,T]\times\R^d\to\R)$ be such that $\varphi-w$ has a local maximum of 0 at $(t_0,x_0)$ with $t_0<T$. Then $\tilde \varphi(t,x):=\varphi (t,x)-\theta e^{-\lambda t}(1+((|x|-K_\Gamma)^+)^{2\gamma+2})\in C^{1,2}([0,T]\times\R^d\to\R)$ and $\tilde \varphi-v$ has a local maximum of 0 at $(t_0,x_0)$. Since $v$ is a viscosity supersolution, we have
\begin{align*}
    0&\leq -\partial_t \tilde \varphi(t,x)-\inf_{\alpha\in A}H(t,x,\tilde \varphi(t,x),D\tilde \varphi(t,x),D^2\tilde \varphi(t,x),\alpha)
    \\
    &=-\partial_t(\varphi(t,x)-\theta e^{-\lambda t}(1+((|x|-K_\Gamma)^+)^{2\gamma+2}))-\inf_{\alpha\in A}\big\{\mcL^\alpha (\varphi(t,x)-\theta e^{-\lambda t}(1+((|x|-K_\Gamma)^+)^{2\gamma+2}))
    \\
    &\quad+f(t,x,\varphi(t,x)-\theta e^{-\lambda t}(1+((|x|-K_\Gamma)^+)^{2\gamma+2}),\sigma^\top(t,x)\nabla_x (\varphi(t,x)-\theta e^{-\lambda t}(1+((|x|-K_\Gamma)^+)^{2\gamma+2})),\alpha)\big\}
    \\
    &\leq -\partial_t\varphi(t,x)-\inf_{\alpha\in A}\big\{\mcL^\alpha \varphi(t,x)+f(t,x,\varphi(t,x),\sigma^\top(t,x)\nabla_x \varphi(t,x),\alpha)\big\}
    \\
    &\quad-\theta \lambda e^{-\lambda t}(1+((|x|-K_\Gamma)^+)^{2\gamma+2})+\sup_{\alpha\in A}\mcL^\alpha\{\theta e^{-\lambda t}(1+((|x|-K_\Gamma)^+)^{2\gamma+2})\}
    \\
    &\quad+k_f\theta e^{-\lambda t}(1+((|x|-K_\Gamma)^+)^{2\gamma+2}+C(1+|x|) ((|x|-K_\Gamma)^+)^{2\gamma+1})
\end{align*}
Consequently,
\begin{align*}
&-\partial_t\varphi(t,x)-\inf_{\alpha\in A}\mcL^\alpha \{\varphi(t,x)+f(t,x,\varphi(t,x),\sigma^\top(t,x)\nabla_x \varphi(t,x),\alpha)\}
\\
%&\geq \theta \lambda e^{-\lambda t}((|x|-K_\Gamma)^+)^{2\gamma+2}-\theta\mcL e^{-\lambda t}((|x|-K_\Gamma)^+)^{2\gamma+2}
%\\
%&\quad f(t,x,\varphi(t,x)-\theta e^{-\lambda t}((|x|-K_\Gamma)^+)^{2\gamma+2},\sigma^\top(t,x)\nabla_x (\varphi(t,x)-\theta e^{-\lambda t}((|x|-K_\Gamma)^+)^{2\gamma+2}))
%\\
%&\quad-f(t,x,\varphi(t,x),\sigma^\top(t,x)\nabla_x \varphi(t,x))
%\\
&\geq \theta e^{-\lambda t}\big(\lambda(1+((|x|-K_\Gamma)^+)^{2\gamma+2})-C(1+|x|)((|x|-K_\Gamma)^+)^{2\gamma+1}- C(1+|x|)^2 e^{-\lambda t}((|x|-K_\Gamma)^+)^{2\gamma}
\\
&\quad -k_f (1+((|x|-K_\Gamma)^+)^{2\gamma+2}+C(1+|x|) ((|x|-K_\Gamma)^+)^{2\gamma+1})\big),
\end{align*}
where the right hand side is non-negative for all $\theta> 0$ and all $\lambda>\lambda_0$ for some $\lambda_0>0$.\qed\\

We have the following results the proof of which we omit since it is classical:
\begin{lem}\label{lem:integ-factor}
A locally bounded function $v:[0,T]\times \R^d\to\R$ is a viscosity supersolution (resp. subsolution) to \eqref{ekv:var-ineq} if and only if for every $\lambda\in\R$, $\tilde v(t,x):=e^{\lambda t}v(t,x)$ is a viscosity supersolution (resp. subsolution) to
\begin{align}\label{ekv:var-ineq-if}
\begin{cases}
  \min\big\{\tilde v(t,x)-\sup_{b\in U}\{\tilde v(t,\Gamma(t,x,b))-e^{\lambda t}\ell(t,x,b)\},-\tilde v_t(t,x)+\lambda \tilde v(t,x)-\inf_{\alpha\in A}\{\mcL^\alpha \tilde v(t,x)\\
  +e^{\lambda t}f(t,x,e^{-\lambda t}\tilde v(t,x),e^{-\lambda t}\sigma^\top(t,x)\nabla_x \tilde v(t,x),\alpha)\}\big\}=0,\quad\forall (t,x)\in[0,T)\times \R^d \\
  \tilde v(T,x)=e^{\lambda T}\psi(x).
\end{cases}
\end{align}
\end{lem}
\begin{rem}
Here, it is important to note that $\tilde \ell(t,x):=e^{\lambda t}\ell(t,x)$, $\tilde f(t,x,y,z,\alpha):=-\lambda y\\+e^{\lambda t}f(t,x,e^{-\lambda t}y,e^{-\lambda t}z,\alpha)$ and $\tilde \psi(x):=e^{\lambda T}\psi(x)$ satisfy Assumption~\ref{ass:on-coeff}. In particular, this implies that Lemma~\ref{lem:is-super} holds for supersolutions to \eqref{ekv:var-ineq-if} as well.
\end{rem}

We have the following comparison result for viscosity solutions in $\Pi_{pg}$:

\begin{prop}
Let $v$ (resp. $u$) be a supersolution (resp. subsolution) to \eqref{ekv:var-ineq}. If $u,v\in \Pi_{pg}$, then $u\leq v$.
\end{prop}

\noindent \emph{Proof.} First, we note that we only need to show that the statement holds for solutions to~\eqref{ekv:var-ineq-if}. We thus assume that $v$ (resp.~$u$) is a viscosity supersolution (resp.~subsolution) to \eqref{ekv:var-ineq-if}. %Moreover, we may without loss of generality assume that $v$ is l.s.c.~and $u$ is u.s.c.

It is sufficient to show that
\begin{align*}
w(t,x)&=w^{\theta,\lambda}(t,x):=v(t,x)-\theta e^{-\lambda t}(1+((|x|-K_\Gamma)^+)^{2\gamma+2})
\\
&\geq u(t,x)
\end{align*}
for all $(t,x)\in[0,T]\times\R^d$ and any $\theta>0$. Then the result follows by taking the limit $\theta\to 0$. Moreover, we know from Lemma~\ref{lem:is-super} that there is a $\lambda_0>0$ such that $w$ is a supersolution to \eqref{ekv:var-ineq-if} for each $\lambda\geq\lambda_0$ and $\theta>0$.

By assumption, $u,v\in \Pi_{pg}$, which implies that there are $C>0$ and $\gamma>0$ such that
\begin{align*}
|v(t,x)|+|u(t,x)|\leq C(1+|x|^{2\gamma}).
\end{align*}
Hence, for each $\lambda,\theta>0$ there is a $R\geq K_\Gamma$ such that
\begin{align*}
w(t,x)>u(t,x),\quad\forall (t,x)\in[0,T]\times\R^d,\:|x|>R.
\end{align*}
We search for a contradiction and assume that there is a $(t_0,x_0)\in [0,T]\times \R^d$ such that $v(t_0,x_0)>w(t_0,x_0)$. Then there is a point $(\bar t,\bar x)\in[0,T)\times B_R$ (the open unit ball of radius $R$ centered at 0) such that
\begin{align*}
\max_{(t,x)\in[0,T]\times\R^d}(u(t,x)-w(t,x))&=\max_{(t,x)\in[0,T)\times B_R}(u(t,x)-w(t,x))
\\
&=u(\bar t,\bar x)-w(\bar t,\bar x)=\eta>0.
\end{align*}
We first show that there is at least one point $(t^*,x^*)\in[0,T)\times B_R$ such that
\begin{enumerate}[a)]
  \item $u(t^*,x^*)-w(t^*,x^*)=\eta$ and
  \item $u(t^*,x^*)>\sup_{b\in U}\{u(t^*,\Gamma(t^*,x^*,b))-\tilde\ell(t^*,b)\}$.
\end{enumerate}
We again argue by contradiction and assume that $u(t,x)=\sup_{b\in U}\{u(t,\Gamma(t,x,b))-\tilde\ell(t,b)\}$ for all $(t,x)\in A:=\{(s,y)\in[0,T]\times\R^d: u(s,y)-w(s,y)=\eta\}$. Indeed, as $u$ is u.s.c. and $\Gamma$ is continuous, there is a $b_1$ such that
\begin{align}\label{ekv:equal}
u(\bar t,\bar x)=\sup_{b\in U}\{u(\bar t,\Gamma(\bar t,\bar x,b))-\tilde\ell(\bar t,b)\}=u(\bar t,\Gamma(\bar t,\bar x,b_1))-\tilde\ell(\bar t,b_1).
\end{align}
Now, set $x_1=\Gamma(\bar t,\bar x,b_1)$ and note that since
 \begin{align*}
|\Gamma(t,x,b)|\leq R,\quad \forall(t,x,b)\in [0,T]\times \bar B_R\times U,
\end{align*}
it follows that $x_1\in \bar B_R$. Moreover, as $w$ is a supersolution it satisfies
\begin{align*}
  w(\bar t,\bar x)- (w(\bar t,\Gamma(\bar t,\bar x,b_1))-\tilde\ell(t,\bar x,b_1))\geq 0
\end{align*}
or
\begin{align*}
  - w(\bar t,x_1))\geq -w(\bar t,\bar x)-\tilde\ell(t,\bar x,b_1)
\end{align*}
and we conclude from \eqref{ekv:equal} that
\begin{align*}
  u(\bar t,x_1)- w(\bar t,x_1)&\geq u(\bar t,\bar x)+\tilde\ell(\bar t,\bar x,b_1)-(w(\bar t,\bar x)+\tilde\ell(\bar t,\bar x,b_1))
  \\
  &=u(\bar t,\bar x)-w(\bar t,\bar x)=\eta.
\end{align*}
Hence, $(\bar t,x_1)\in A$ and by our assumption it follows that there is a $b_2\in U$ such that
\begin{align*}%\label{ekv:equal}
u(\bar t,x_1)=u(\bar t,\Gamma(\bar t,x_1,b_2))-\tilde\ell(\bar t,b_2)
\end{align*}
and a corresponding $x_2:=\Gamma(\bar t,x_1,b_2)\in B_R$. Now, this process can be repeated indefinitely to find a sequence $(x_j,b_j)_{j\geq 1}$ in $B_R\times U$ such that for any $l\geq 0$ we have
\begin{align*}
  u(\bar t,\bar x)=u(\bar t,x_l)-\sum_{j=1}^{l}\tilde\ell(\bar t,x_{j-1},b_j),
\end{align*}
with $x_0:=\bar x$. Now, as $\tilde\ell\geq (1\wedge e^{\lambda T})\delta>0$ we get a contradiction by letting $l\to\infty$ while noting that $|u(t,x)|$ is bounded on $[0,T]\times \bar B_R$. We can thus pick a $(t^*,x^*)\in [0,T)\times B_R$ such that \emph{a)} and \emph{b)} above holds.

The remainder of the proof is similar to the proof of Proposition 4.1 in~\cite{Morlais13}. We assume without loss of generality that $\gamma\geq 2$ and define
\begin{align*}
\Phi_n(t,x,y):=u(t,x)-w(t,x)-\varphi_n(t,x,y),
\end{align*}
where% for some $\bar\theta\in (0,1]$
\begin{align*}
  \varphi_n(t,x,y):=\frac{n}{2}|x-y|^{2\gamma}+|x-x^*|^2+|y-x^*|^2+(t-t^*)^2.
\end{align*}
Since $u$ is u.s.c.~and $w$ is l.s.c.~there is a $(t_n,x_n,y_n)\in[0,T]\times \bar B_R\times \bar B_R$ (with $\bar B_R$ the closure of $B_R$) such that
\begin{align*}
  \Phi_n(t_n,x_n,y_n)=\max_{(t,x,y)\in [0,T]\times \bar B_R\times \bar B_R}\Phi_n(t,x,y).
\end{align*}
Now, the inequality $2\Phi_n(t_n,x_n,y_n)\geq \Phi_n(t_n,x_n,x_n)+\Phi_n(t_n,y_n,y_n)$ gives
\begin{align*}
n|x_n-y_n|^{2\gamma}\leq u(t_n,x_n)-u(t_n,y_n)+w(t_n,x_n)-w(t_n,y_n).
\end{align*}
Consequently, $n|x_n-y_n|^{2\gamma}$ is bounded (since $u$ and $w$ are bounded on $[0,T]\times\bar B_R\times\bar B_R$) and $|x_n-y_n|\to 0$ as $n\to\infty$. We can, thus, extract subsequences $n_l$ such that $(t_{n_l},x_{n_l},y_{n_l})\to (\tilde t,\tilde x,\tilde x)$ as $l\to\infty$. Since
\begin{align*}
u(t^*,x^*)-w(t^*,x^*)\leq \Phi_n(t_n,x_n,y_n)\leq u(t_n,x_n)-w(t_n,y_n),
\end{align*}
it follows that
\begin{align*}
u(t^*,x^*)-w(t^*,x^*)&\leq \limsup_{l\to\infty} \{u(t_{n_l},x_{n_l})-w(t_{n_l},y_{n_l})\}
\\
&\leq u(\tilde t,\tilde x)-w(\tilde t,\tilde x)
\end{align*}
and as the righthand side is dominated by $u(t^*,x^*)-w(t^*,x^*)$ we conclude that
\begin{align*}
  u(\tilde t,\tilde x)-w(\tilde t,\tilde x)=u(t^*,x^*)-w(t^*,x^*).
\end{align*}
In particular, this gives that $\lim_{l\to\infty}\Phi_n(t_{n_l},x_{n_l},y_{n_l})=u(\tilde t,\tilde x)-w(\tilde t,\tilde x)$ which implies that
\begin{align*}
  \limsup_{l\to\infty} n_l|x_{n_l}-y_{n_l}|^{2\gamma}= 0
\end{align*}
and
\begin{align*}
  (t_{n_l},x_{n_l},y_{n_l})\to (t^*,x^*,x^*).
\end{align*}
We can extract a subsequence $(\tilde n_l)_{l\geq 0}$ of $(n_l)_{l\geq 0}$ such that $t_{\tilde n_l}<T$, $|x_{\tilde n_l}|<R$ and
\begin{align*}
  u(t_{\tilde n_l},x_{\tilde n_l})-w(t_{\tilde n_l},x_{\tilde n_l})\geq \frac{\eta}{2}.
\end{align*}
Moreover, since $\sup_{b\in U}\{u(t,\Gamma(t,x,b))-\tilde\ell(t,b)\}$ is u.s.c.~(see Lemma~\ref{lem:mcM-monotone}) and $u(t_{\tilde n_l},x_{\tilde n_l})\to u(t^*,x^*)$ there is an $l_0\geq 0$ such that
\begin{align*}
  u(t_{\tilde n_l},x_{\tilde n_l})-\sup_{b\in U}\{u(t_{\tilde n_l},\Gamma(t_{\tilde n_l},x_{\tilde n_l},b))-\tilde\ell(t_{\tilde n_l},b)\}>0,
\end{align*}
for all $l\geq l_0$. To simplify notation we will, from now on, denote $(\tilde n_l)_{l\geq l_0}$ simply by $n$.\\

By Theorem 8.3 of~\cite{UsersGuide} there are $(p^u_n,q^u_n,M^u_n)\in \bar J^{2,+}u(t_n,x_n)$ and $(p^w_n,q^w_n,M^w_n)\in \bar J^{2,+}w(t_n,y_n)$
such that
\begin{align*}
\begin{cases}
  p^u_n-p^w_n=\partial_t\varphi_n(t_n,x_n,y_n)=2(t_n-t^*)
  \\
  q^u_n=D_x\varphi_n(t_n,x_n,y_n)=n\gamma(x-y)|x-y|^{2\gamma-2}+2(x-x^*)
  \\
  q^w_n=-D_y\varphi_n(t_n,x_n,y_n)=n\gamma(x-y)|x-y|^{2\gamma-2}-2(y-x^*)
\end{cases}
\end{align*}
and for every $\epsilon>0$,
\begin{align*}
  \left[\begin{array}{cc} M^n_x & 0 \\ 0 & -M^n_y\end{array}\right]\leq B_n(t_n,x_n,y_n)+\epsilon B_n^2(t_n,x_n,y_n),
\end{align*}
where $B_n(t_n,x_n,y_n):=D^2_{(x,y)}\varphi_n(t_n,x_n,y_n)$. Now, we have
\begin{align*}
  D^2_{(x,y)}\varphi_n(t,x,y)=\left[\begin{array}{cc} D_x^2\varphi_n(t,x,y) & D^2_{yx}\varphi_n(t,x,y) \\ D^2_{xy}\varphi_n(t,x,y) & D_y^2\varphi_n(t,x,y)\end{array}\right]
  = \left[\begin{array}{cc} n\xi(x,y)+2 I & -n\xi(x,y) \\ -n\xi(x,y) & n\xi(x,y) +2 I \end{array}\right]
\end{align*}
where $I$ is the identity-matrix of suitable dimension and
\begin{align*}
  \xi(x,y):=\gamma|x-y|^{2\gamma-4}\{|x-y|^2I+2(\gamma-1)(x-y)(x-y)^\top\}.
\end{align*}
In particular, since $x_n$ and $y_n$ are bounded, choosing $\epsilon:=\frac{1}{n}$ gives that
\begin{align}\label{ekv:mat-bound}
  \tilde B_n:=B_n(t_n,x_n,y_n)+\epsilon B_n^2(t_n,x_n,y_n)\leq Cn|x_n-y_n|^{2\gamma-2}\left[\begin{array}{cc} I & -I \\ -I & I \end{array}\right]+C I.
\end{align}
By the definition of viscosity supersolutions and subsolutions we have that
\begin{align*}
&-p^u_n+\lambda u(t_n,x_n)-a^\top(t_n,x_n,\alpha)q^u_n-\frac{1}{2}\trace [\sigma^\top(t_n,x_n,\alpha)M^u_n\sigma(t_n,x_n,\alpha)]
\\
&-e^{\lambda t_n}f(t_n,x_n,e^{-\lambda t_n}u(t_n,x_n),e^{-\lambda t_n}\sigma^\top(t_n,x_n)q^u_n,\alpha)\big\}\leq 0
\end{align*}
for all $\alpha\in A$ and
\begin{align*}
&-p^w_n+\lambda w(t_n,y_n)-\inf_{\alpha\in A}\big\{a^\top(t_n,y_n,\alpha)q^w_n+\frac{1}{2}\trace [\sigma^\top(t_n,y_n,\alpha)M^w_n\sigma(t_n,y_n,\alpha)]
\\
&+e^{\lambda t_n}f(t_n,y_n,e^{-\lambda t_n}w(t_n,y_n),e^{-\lambda t_n}\sigma^\top(t_n,x_n)q^w_n,\alpha)\big\}\geq 0.
\end{align*}
Combined, this gives that
\begin{align*}
\lambda (u(t_n,x_n)-w(t_n,y_n))&\leq \sup_{\alpha\in A}\big\{p^u_n+a^\top(t_n,x_n,\alpha)q^u_n+\frac{1}{2}\trace [\sigma^\top(t_n,x_n,\alpha)M^u_n\sigma(t_n,x_n,\alpha)]
\\
&+e^{\lambda t_n}f(t_n,x_n,e^{-\lambda t_n}u(t_n,x_n),e^{-\lambda t_n}\sigma^\top(t_n,x_n)q^u_n,\alpha)
\\
&-p^w_n-a^\top(t_n,y_n,\alpha)q^w_n-\frac{1}{2}\trace [\sigma^\top(t_n,y_n,\alpha)M^w_n\sigma(t_n,y_n,\alpha)]
\\
&-e^{\lambda t_n}f(t_n,y_n,e^{-\lambda t_n}w(t_n,y_n),e^{-\lambda t_n}\sigma^\top(t_n,x_n)q^w_n,\alpha)\big\}
\end{align*}
Collecting terms we have that
\begin{align*}
p^u_n-p^w_n&=2(t_n-t^*)
\end{align*}
and since $a$ is Lipschitz continuous in $x$ and bounded on $\bar B_R$, we have
\begin{align*}
  a^\top(t_n,x_n,\alpha)q^u_n-a^\top(t_n,y_n,\alpha)q^w_n&\leq  (a^\top(t_n,x_n,\alpha)-a^\top(t_n,y_n,\alpha))n\gamma(x_n-y_n)|x_n-y_n|^{2\gamma-2}
  \\
  &\quad+C(|x_n-x^*|+|y_n-x^*|)
  \\
  &\leq C(n|x_n-y_n|^{2\gamma}+|x_n-x^*|+|y_n-x^*|),
\end{align*}
where the right-hand side tends to 0 as $n\to\infty$. Let $s_x$ denote the $i^{\rm th}$ column of $\sigma(t_n,x_n,\alpha)$ and let $s_y$ denote the $i^{\rm th}$ column of $\sigma(t_n,y_n,\alpha)$ then by the Lipschitz continuity of $\sigma$ and \eqref{ekv:mat-bound}, we have
\begin{align*}
  s_x^\top M^u_n s_x-s_y^\top M^w_n s_y&=\left[\begin{array}{cc} s_x^\top  & s_y^\top \end{array}\right]\left[\begin{array}{cc} M^u_n  & 0 \\ 0 &-M^w_n\end{array}\right]\left[\begin{array}{c} s_x \\ s_y \end{array}\right]
  \\
  &\leq \left[\begin{array}{cc} s_x^\top  & s_y^\top \end{array}\right]\tilde B_n\left[\begin{array}{c} s_x \\ s_y \end{array}\right]
  \\
  &\leq C(n|x_n-y_n|^{2\gamma}+|x_n-y_n|)
\end{align*}
and we conclude that
\begin{align*}
  \limsup_{n\to\infty}\sup_{\alpha\in A}\frac{1}{2}\trace [\sigma^\top(t_n,x_n,\alpha)M^u_n\sigma(t_n,x_n,\alpha)-\sigma^\top(t_n,y_n,\alpha)M^w_n\sigma(t_n,y_n,\alpha)]\leq 0.
\end{align*}
Finally, we have for some $C_R>0$ that
\begin{align*}
  &e^{\lambda t_n}f(t_n,x_n,e^{-\lambda t_n}u(t_n,x_n),e^{-\lambda t_n}\sigma^\top(t_n,x_n)q^u_n,\alpha)-e^{\lambda t_n}f(t_n,y_n,e^{-\lambda t_n}w(t_n,y_n),e^{-\lambda t_n}\sigma^\top(t_n,x_n)q^w_n,\alpha)
  \\
  &\leq k_f(u(t_n,x_n)-w(t_n,y_n))+C_R(|x_n-y_n|+|\sigma^\top(t_n,x_n,\alpha)q^u_n-\sigma^\top(t_n,x_n,\alpha)q^w_n|).
\end{align*}
Repeating the above argument we get that the upper limit of the right-hand side when $n\to\infty$ is bounded by $k_f(u(t_n,x_n)-w(t_n,y_n))$. Put together, this gives that
\begin{align*}
(\lambda-k_f) \limsup_{n\to\infty}(u(t_n,x_n)-w(t_n,y_n))&\leq 0
\end{align*}
a contradiction since $\lambda\in\R$ was arbitrary.\qed\\

\bibliographystyle{plain}
\bibliography{ImpVCGames_ref}

\begin{thebibliography}{10}

\bibitem{Azimzadeh19}
P.~Azimzadeh.
\newblock A zero-sum stochastic differential game with impulses,precommitment,
  and unrestricted cost functions.
\newblock {\em Appl Math Optim}, 79:483--514, 2019.

\bibitem{BayraktarRobust}
E.~Bayraktar, A.~Cosso, and H.~Pham.
\newblock Robust feedback switching control: dynamic programming and viscosity
  solutions.
\newblock {\em SIAM J. Control Optim.}, 54(5):2594--2628, 2016.

\bibitem{BensLionsImpulse}
A.~Bensoussan and J.L. Lions.
\newblock {\em Impulse Control and Quasivariational inequalities}.
\newblock Gauthier-Villars, Montrouge, France, 1984.

\bibitem{Buckdahn08}
R.~Buckdahn and J.~Li.
\newblock Stochastic differential games and viscosity solutions of
  hamilton-jacobi-bellman-isaacs equations.
\newblock {\em SIAM J. Control Optim}, 47(1):444--475, 2008.

\bibitem{Cosso13}
A.~Cosso.
\newblock Stochastic differential games involving impulse controls and
  double-obstacle quasi-variational inequalities.
\newblock {\em SIAM J. Control Optim.}, 3(51):2102--2131, 2013.

\bibitem{UsersGuide}
M.~G. Crandall, H.~Ishii, and P.~L. Lions.
\newblock User's guide to viscosity solutions of second order partial
  differential equations.
\newblock {\em Bulletin of the American Mathematical Society}, 27(1):1--67,
  1992.

\bibitem{BollanSWG1}
B.~Djehiche, S.~Hamad\'ene, and M.~Morlais.
\newblock Viscosity solutions of systems of variational inequalities with
  interconnected bilateral obstacles.
\newblock {\em Funkcialaj Ekvacioj}, 58(1):135--175, 2015.

\bibitem{DjehicheSWG2}
B.~Djehiche, S.~Hamad\'ene, M.-A. Morlais, and X.~Zhao.
\newblock On the equality of solutions of max-min and min-max systems of
  variational inequalities with interconnected bilateral obstacles.
\newblock {\em J. Math. Anal. Appl.}, 452:148--175, 2017.

\bibitem{ElKaroui2}
N.~{El Karoui}, S.~Peng, and M.~C. Quenez.
\newblock Backward stochastic differential equationsin finance.
\newblock {\em Math. Finance}, 7(1):1--71, 1997.

\bibitem{ElliotKalton72}
R.~J. Elliott and N.~J. Kalton.
\newblock {\em The existence of value in differential games}.
\newblock Number 126. Memoirs of the American Mathematical Society, Providence,
  Rhode Island, 1972.

\bibitem{Evans84}
L.~C. Evans and P.~E. Souganidis.
\newblock Differential games and representation formulasfor solutions of
  hamilton-jacobi-isaacs equations.
\newblock {\em Indiana Univ. Math. J.}, 33:773--797, 1984.

\bibitem{FlemSoug89}
W.~H. Flemming and P.~E. Souganidis.
\newblock On the existence of value functions of two-player, zero-sum
  stochastic differential games.
\newblock {\em Indiana Univ. Math. J.}, 38:293--314, 1989.

\bibitem{Morlais13}
S.~Hamad{\'e}ne and M.~A. Morlais.
\newblock Viscosity solutions of systems of pdes with interconnected obstacles
  and switching problem.
\newblock {\em Appl Math Optim.}, 67:163--196, 2013.

\bibitem{HamZhang}
S.~Hamad\'ene and J.~Zhang.
\newblock Switching problem and related system of reflected backward {SDEs}.
\newblock {\em Stochastic Process. Appl.}, 120(4):403--426, 2010.

\bibitem{HuTang}
Y.~Hu and S.~Tang.
\newblock Multi-dimensional {BSDE} with oblique reflection and optimal
  switching.
\newblock {\em Prob. Theory and Related Fields}, 147(1-2):89--121, 2008.

\bibitem{Isaacs65}
R.~Isaacs.
\newblock {\em Differential games. A mathematical theory with applications to
  warfare andpursuit, control and optimization}.
\newblock John Wiley \& Sons, Inc., New York-London-Sydney, 1965.

\bibitem{LiPeng09}
J.~Li and S.~Peng.
\newblock Stochastic optimization theory of backward stochastic
  differential equations with jumps and viscosity solutions
  ofhamilton-jacobi-bellman equations.
\newblock {\em Nonlinear Analysis}, 70:1776--1796, 2009.

\bibitem{Li14}
J.~Li and Q.~Wei.
\newblock Optimal control problems of fully coupled fbsdes and viscosity
  solutions of hamilton-jacobi-bellman equations.
\newblock {\em SIAM J. Control Optim}, 52(3):1622--1662, 2014.

\bibitem{rbsde_impulse}
M.~Perninge.
\newblock Finite horizon robust impulse control in a non-markovian framework
  and related systems of reflected bsdes.
\newblock {\em arXiv:2103.16272}, 2021.

\bibitem{Protter}
P.~Protter.
\newblock {\em Stochastic Integration and Differential Equations}.
\newblock Springer, Berlin, 2nd edition, 2004.

\bibitem{TangHouSWgame}
S.~Tang and {Sh}. Hou.
\newblock Switching games of stochastic differential systems.
\newblock {\em SIAM J. Control Optim.}, 46(3):900--929, 2007.

\bibitem{FZhang11}
F.~Zhang.
\newblock Stochastic differential games involving impulse controls.
\newblock {\em ESAIM Control Optim. Calc. Var.}, 17(3):749--760, 2011.

\bibitem{LZhang21}
L.~Zhang.
\newblock A bsde approach to stochastic differential games involvingimpulse
  controls and hjbi equation.
\newblock {\em J Syst Sci Complex}, 2021.

\end{thebibliography}
\end{document}